\documentclass[12pt,a4paper]{amsart}
\usepackage{xcolor}
\usepackage[draft]{todonotes} 
\usepackage[english]{babel}
\usepackage{geometry}

\hoffset=- 1cm \voffset=1cm \topmargin=-0.5in \textheight=23cm
\usepackage[T1]{fontenc}
\usepackage[latin1]{inputenc}
\usepackage{amsmath,amssymb,mathrsfs,amsthm, tikz-cd,mathrsfs}
\usepackage{hyperref}
\usepackage{enumitem}
\usepackage{makecell}
\usepackage{datetime}
\usepackage{xcolor}
\usepackage{fancyhdr}

\usepackage{DLdef1}
\newcommand {\fbrj}{\mathfrak{brj}}

\begin{document}

\title{Double extensions of restricted Lie (super)algebras}

\author{Sa\"id Benayadi${}^a$}
\address{${}^a$Laboratoire de Math\'ematiques IECL UMR CNRS 7502,
Universit\'e de Lorraine, 3 rue Augustin Fresnel, BP 45112,
F-57073 Metz Cedex 03, FRANCE.}

\email{said.benayadi@univ-lorraine.fr }

\author{Sofiane Bouarroudj${}^{b,*}$}

 \address{${}^b$Division of Science and Mathematics, New York University Abu Dhabi, Po Box 129188, Abu Dhabi, United Arab Emirates.\\${}^{*}$ The corresponding author}
            
\email{sofiane.bouarroudj@nyu.edu}

\author{Mounir Hajli${}^c$}

\address{${}^c$ School of Mathematical Sciences, Shanghai Jiao Tong University, P.R. China}

\email{hajli@sjtu.edu.cn}






\begin{abstract}

A \textit{double extension} ($\mathscr{D}$-extension) of a~Lie (super)algebra $\fa$ with a~non-degenerate invariant symmetric bilinear form $\mathscr{B}$, briefly: a~\textit{NIS}-(super)algebra, is an enlargement of $\fa$ by means of a~central extension and a~derivation; the affine Kac-Moody algebras are the best known examples of double extensions of loops algebras.  

Let $\fa$ be a~restricted Lie (super)algebra with a~NIS $\mathscr{B}$. Suppose $\fa$ has a~restricted derivation $\mathscr{D}$ such that  $\mathscr{B}$ is $\mathscr{D}$-invariant. We show that the double extension of $\fa$ constructed by means of $\mathscr{B}$ and $\mathscr{D}$ is restricted. We show that, the other way round, any restricted NIS-(super)algebra with non-trivial center  can be obtained as a~$\mathscr{D}$-extension of another restricted NIS-(super)algebra subject to an extra condition on the central element.

We give new examples of $\mathscr{D}$-extensions of restricted Lie (super)algebras, and pre-Lie superalgebras indigenous to characteristic 3.
\end{abstract}


\maketitle

\thispagestyle{empty}
\setcounter{tocdepth}{2}
\tableofcontents

\section{Introduction} \label{SecDef}

\subsection{NIS-superalgebras}\label{Nises}
Let $\mathscr{B}$ be a~bilinear form on $\fa$. Consider the \textit{upsetting} of bilinear forms 
$u\colon\text{Bil} (V, W)\longmapsto\text{Bil}(W, V)$, see \cite[Ch.1]{SoS}, given by the formula 
\begin{equation}\label{susyB}
u(B)(w, v)=(-1)^{p(v)p(w)}B(v, w)\text{~~ for any $v \in V$ and $w\in W$.}
\end{equation}

If $V=W$, we say that
$B$ is \textit{symmetric} (resp. \textit{anti-symmetric}) if  $u(B)=B$ (resp. $u(B)=-B$).

Following \cite{BKLS}, we call the Lie superalgebra $\fa$ a~\textit{NIS-superalgebra} (sometimes used to be called \textit{quadratic} in the literature) if it has a~non-degenerate, invariant, and symmetric bilinear form $\mathscr{B}$. We denote such a~superalgebra by $(\fa, \mathscr{B})$.

For a~list of a~wide class of simple modular NIS-superalgebras, see \cite{BKLS}.

A NIS-superalgebra $(\fa, \mathscr{B})$ is said to be \textit{reducible} if it can be decomposed into a~direct sum of mutually orthogonal ideals, namely $\fa=\oplus I_i$; otherwise, it is called \textit{irreducible}.

\subsection{Double extensions in general} Let $\fa$ be a~NIS-superalgebra (in particular, a~Lie algebra) 
defined over a~field $\mathbb{K}$ of positive characteristic $p$. The notion of a~double extension of the Lie superalgebra $\fa$, called \textit{$\mathscr{D}$-extension} in \cite{BeBou}, was introduced by Medina and Revoy \cite{MR} in the case of Lie algebras over $\mathbb{R}$. This notion has been \textit{superized} and studied in a~series of papers \cite{ABB, ABBQ, BBB, BB, B, B2, BeBou};  for a~succinct summary, see \cite{BLS}.

The double extension of $\fa$, denoted by $\fg$,  simultaneously involves three ingredients:
\begin{itemize}
\item a~derivation $\mathscr{D}$,
\item a~$\mathscr{D}$-invariant NIS $\mathscr{B}_\fa$ on $\fa$,
\item the central extension $\fa_x $ of $\fa$ is given by the cocycle $(a,b)\mapsto\mathscr{B}_\fa(\mathscr{D}(a),b)x$ for any $a,b\in\fa$ with the center spanned by the vector $x$.
\end{itemize}

It turns out that for $\fg$ to be non-trivial, i.e., not isomorphic to the direct sum of $\fa$ and 2-dimensional ideal $\Bbb Kx\oplus \Bbb K \mathscr{D}$, the central extension must be non-trivial and the derivation $\mathscr{D}$ outer, moreover, if $\mathscr{D}$ is odd, the following conditions must be met:
\begin{equation}\label{2cond}
\begin{array}{l}

\text{ $\mathscr{D}^2=\ad_b$ and $\mathscr{D}(b)=0$  for some $b\in \fa_\ev$.}
\end{array}
\end{equation}



\textit{Modular} Lie (super)algebras, i.e., Lie (super)algebras defined over  the field of characteristic $p>0$, can be double extended. The inductive description \`a la Medina and Revoy becomes, however, the most challenging part, because Lie's theorem and the Levi decomposition do not hold if $p>0$. The formulations of the inductive description differ if $p>0$ from their counterparts for $p=0$, see Theorems \ref{Rec1}, \ref{DeBeRec},~\ref{DoBoRec}.

Favre and Santharoubane \cite{FS} introduced an important ingredient in the study of NIS-algebras $\fg$ and $\tilde \fg$ -- double extensions of the same Lie algebra~$\fa$: they suggested to consider them up to an isomorphism $\pi : \fg \rightarrow \tilde \fg$ (Favre and Santharoubane called $\pi$ an isometry) such that
\[
\mathscr{B}_{\tilde \fg} (\pi(f ), \pi(g)) = \mathscr{B}_\fg(f, g) \text{ for any } f, g \in \fg.
\] The equivalence of double extensions up to isomorphisms turns out to be a~very important and useful notion, as demonstrated in \cite{BeBou, BE, BDRS, BLS} by several new examples. In \S \ref{subsubsec:SIsom}, we generalize this notion to the case of Lie superalgebras for any $p>2$, completing the results of \cite{BeBou} where $p=2$ was considered. 

Observe that for Lie algebras over fields of characteristic $p > 0$, and for Lie superalgebras over fields of any characteristic, the Killing form might be degenerate, cf. \cite{BGP, SF, GP, BKLS}. 

Simple NIS-superalgebras having outer derivations are  
abundant in positive characteristic (see \cite{Dz, Fa, BKLS, B}), and therefore might be double extended if they have non-trivial central extensions; for their list, see \cite{BGLL1}.

In \cite{KO},  an interesting relation of double extensions with pseudo-Riemannian manifolds is established; the two approaches (all the above and \cite{KO}) have an empty intersection.

\subsection{Double extensions of restricted Lie (super)algebras}
As far as we know, the notion of a~\textit{restricted} Lie algebra was introduced by Jacobson \cite{J}. Roughly speaking, one requires the existence of an endomorphism on the modular Lie algebra that resembles the $p$-th power mapping $x\mapsto x^p$ in associative algebras. Lie algebras associated with algebraic groups over fields of positive characteristic are  
restricted, and this class resembles the characteristic $0$ case, see \cite{SF, S}. Superization of the notion of restrictedness has been studied by several authors, see~\cite{P, BKLLS, BGL, Fa2} and especially \cite{BLLSq} where new phenomena were observed. Here we consider only the classical restrictedness, see \S \ref{defres}.

Let $(\fa,\mathscr{B}_\fa)$ be a~restricted NIS-(super)algebra with a~derivation $\mathscr{D}$ such that 
$\mathscr{B}_\fa$  is $\mathscr{D}$-invariant, namely, 
\begin{equation}\label{SD-inv}
\begin{array}{l}
\mathscr{B}_\fa(\mathscr{D}(a),b)+(-1)^{p(\mathscr{D})p(a)}\mathscr{B}_\fa(a,\mathscr{D}(b))=0 \text{ for all $a,b\in \fa$ and $p\neq 2$}; \\[2mm] 
\text{for $p=2$, we additionally require}\\
\text{$\mathscr{B}_\fa(a,\mathscr{D}(a))=0$ for all $a\in \fa$, whenever $p(\mathscr{D})+p(\mathscr{B}_\fa)=\ev$.}
\end{array}
\end{equation}
In particular, if $\mathscr{D}=\ad_c$, then
\[
\mathscr{B}_\fa([c,a],b)+(-1)^{p(c)p(a)}\mathscr{B}_\fa(a,[c,b])=0\text{~~for any $a,b,c\in \fa$.}
\]
In this paper we answer under what condition
\begin{enumerate}
\item [(i)] the $p$-structure on $\fa$ with NIS $\mathscr{B}_\fa$ can be extended to a~$p$-structure on the double extension $\fg$ of $(\fa, \mathscr{B}_\fa)$. 
\item[(ii)] a~restricted NIS-superalgebra $\fg$ is a~
double extension of a~restricted NIS-(super)al\-geb\-ra~$\fa$.
\end{enumerate}

\underline{Conditions for (i)}: for all characteristics, provided the derivation $\mathscr{D}$ is restricted (see condition (\ref{dvres}) in \S \ref{defres}) and satisfies the following conditions: 

(a) In the case of Lie algebras, $\mathscr{D}$ must satisfy condition (\ref{Conp}), called $p$-\textit{property},  which says that the derivations $\mathscr{D}^p$ and $\gamma \mathscr{D}$, where $\gamma \in \mathbb K$, have to be cohomologous, 
see \S \ref{Rod};  

(b) In the case of Lie superalgebras and $p(\mathscr{D})=\od$, condition \eqref{2cond} is required;
if $p(\mathscr{B})=\ev$, then condition (\ref{StarT}) is required,
see Theorems \ref{DoBe} and \ref{DoBo}.

(c) In the case of Lie superalgebras  and $p(\mathscr{D})=\ev$, the derivation $\mathscr{D}$ must satisfy condition \eqref{SConp}, a~superversion of condition (\ref{Conp}), i.e., have $p$-property, see Theorems  \ref{DoBe} and \ref{DoBo}.

The NIS on $\fg$ is given by eq.~\eqref{Leit}. 

Observe that there is a~large class of Lie algebras whose restricted derivations have $p$-\textit{property} (\ref{Conp}). 
For instance, \textit{nilpotent} restricted Lie algebras, it was proved in \cite{FSW} that there are outer derivations for which $\gamma=0$ and $a_0=0$, see Eq.(\ref{FSW}); for examples of  \textit{simple} restricted Lie algebras with $p$-property, see \S \ref{sec: exa}. 

\underline{Conditions for (ii)}: any restricted NIS-superalgebra $\fg$ can be obtained as a~double extension of a~restricted NIS-superalgebra $\fa$ provided the center of $\fg$ is not trivial, and the orthogonal complement of the central element is a~$p$-ideal, see \S \ref{defres}.

We show that the condition on a~restricted derivation $\mathscr{D}$ of the restricted Lie (super)algebra $\fa$ to  have $p$-\textit{property} \eqref{Conp} or \eqref{SConp} is 
necessary to get a double extension of~$\fa$.

We introduce the notion of \textit{equivalent double extensions} that takes the $p$-structures into account, see Theorems \ref{pIsom1}, \ref{pIsom2}, \ref{pIsom3}. 

If $\fh$ is a restricted Lie algebra, and $p=2$, then the Manin triple $\fg=\fh \oplus \fh^*$ is also restricted, and the  $p$-structure is induced from $\fh$, see \S \ref{Maninp=2}. (Is the same true for $p>2$?) 
\subsection{Examples, see \S~4, summary}
These examples, based on  the list of {\bf known} simple Lie algebras and superalgebras with NIS, see \cite{BKLS}, and the descriptions of their derivations and central extensions, see \cite{BGLL1}, are classified up to an isomorphism. 
The cases admitting double extensions are collected in Table (\ref{Summa}).

\underline{Lie algebras}. Ibraev, see \cite{Ib2}, proved that  only for $\fg=\mathfrak{psl}(3)$, among finite-dimensional Lie algebras $\fg$ with Cartan matrix and their  simple  subquotients, the dimension of the Lie algebra of outer derivations is $>1$. More precisely, $\mathfrak{out}(\mathfrak{psl}(3)) \simeq \mathfrak{psl}(3)$. 

Among known simple restricted $\mathbb{Z}$-graded vectorial Lie algebras, see \cite{BGLLS1}, only $\fsvect^{(1)}(3)$ can be double extended, see \S \ref{vectorial}. Since the vectorial Lie algebras whose shearing vector $\underline N$ has at least one coordinate $>1$ can not be restricted, we assume that $\underline N=(1,\ldots, 1)$. For open cases, see \S~\ref{Secpsl4}. 

The Manin triple $\fg= \fhei(2)\oplus \fhei(2)^*$, although not simple, provides us  with an exceptional (see \eqref{FSW}) example in which the $p$-property (\ref{Conp}) is satisfied with $ \gamma=1$ and $a_0\not =0$, see Eq. (\ref{g1a0}). 

\underline{Lie superalgebras}.  Among Lie superalgebras with Cartan matrix and their simple subquotients, only $\fg(2,3)^{(1)}/\fz$ and $\fg(3,3)^{(1)}/\fz$ can be double extended to Lie superalgebras, see \cite{BGLL1}. Amazingly, the double extensions of $\mathfrak{osp}(1|2)$ and $\mathfrak{brj}(2,3)$ for $p=3$ are \textbf{pre-Lie} superalgebras, not Lie superalgebras, see Appendix. 

\begin{equation}\label{Summa}
\footnotesize
\renewcommand{\arraystretch}{1.4}
\begin{tabular}{|c|c|c|c|c|} \hline
The Lie (super)algebra & Its double extensions&$p$ & $p(\mathscr{B})$\\
\hline
$\fpsl(3)$ & $\fgl(3), \,\widetilde \fgl(3), \,\widehat{\fgl}(3)$  &3& $\ev$\\ \hline

   $\fsvect^{(1)}(3)$ & $\widehat \fsvect^{(1)}(3)$ &3 & $\ev$\\ \hline
   
      $  \fhei(2)\oplus \fhei(2)^*$ & $ \widehat{\fhei(2)\oplus \fhei(2)^*}$ & $2$ &  $\ev$\\ \hline \hline
   
      $\fg(2,3)^{(1)}/\fz$ & $\fg(2,3), \,\widetilde \fg(2,3), \,\widehat{\fg}(2,3)$  &3&  $\ev$\\ \hline

$\fg(3,3)^{(1)}/\fz$ & $\fg(3,3), \,\widetilde \fg(3,3)$  &3&  $\ev$\\ \hline

        $\mathfrak{psq}(n)$ for $n>2$& $\mathfrak{q}(n)$ &$>2$& $\od$\\ \hline
\end{tabular}
\end{equation}

None of the simple vectorial Lie superalgebras considered in \cite{BKLS} can be non-trivially double extended. Certain double extensions of $\mathfrak{psq}(n)$ for $n>2$ and $p=2$ are described in \cite{KLLS}.


\section{Main definitions}
Hereafter, $\mathbb{K}$ is an arbitrary field of characteristic $\mathrm{char}(\mathbb{K})=p$. 
\subsection{Restricted Lie algebras}\label{defres}
Let $\fa$ be a~finite-dimensional modular Lie algebra over~$\mathbb{K}$. For a~comprehensive study of modular 
Lie algebras, see \cite{S, SF}. 

Following \cite{J, SF}, a~mapping $[p]:\fa\rightarrow \fa, \quad a\mapsto a^{[p]}$ is called a~\textit{$p$-structure} of $\fa$ and $\fa$ is said to be {\it restricted}  if 
\begin{equation}\label{RRR}
\begin{array}{l}
\text{$\ad_{a^{[p]}}=(\ad_a)^p$ for all $a \in \fa$;}\\[2pt]

\text{$(\alpha a)^{[p]}=\alpha^p a^{[p]}$ for all $a\in \fa$ and $\alpha \in \mathbb{K}$;}\\[2pt]

\text{$(a+b)^{[p]}=a^{[p]}+b^{[p]}+\sum_{1\leq i\leq p-1}s_i(a,b)$, where the  $s_i(a,b)$ can be obtained from}\\[2mm]
(\ad_{\lambda a+b})^{p-1}(a)=\sum_{1\leq i \leq p-1} is_i(a,b) \lambda^{i-1}.\\
\end{array}
\end{equation}

The following theorem, due to Jacobson, is very useful to us.
\sssbegin{Theorem}[\cite{J}] \label{Jac}
Let $(e_j)_{j\in J}$ be a~basis of $\fa$ such that there are $f_j\in \fa$ satisfying $(\ad_{e_j})^p=\ad_{f_j}$. Then, there exists exactly one $p$-mapping $[p]:\fa\rightarrow \fa$ such that 
\[
e_j^{[p]}=f_j \quad \text{ for all $j\in J$}.
\]
\end{Theorem}
Let $(\fa,[p]_\fa)$ and $(\tilde \fa, [p]_{\tilde \fa})$ be two restricted Lie algebras. A linear map $\pi:\fa \rightarrow \tilde \fa$ is called a~\textit{$p$-homomorphism} if $\pi$ is a~homomorphism of Lie algebras and
\[
\pi(x^{[p]_\fa})=(\pi(x))^{[p]_{\tilde \fa}} \quad \text{for all $x\in \fa$.}
\]

An ideal $I$ of $\fa$ is called a~\textit{$p$-ideal} if $x^{[p]}\in I$ for all $x\in I$. 

For an arbitrary subset $S \subset \fa$, we denote
\[
S^{[p]^i}:=\{x^{[p]^i} \; |\; x\in S  \},
\]
where the expression $[p]^i$ stands for the composition $[p]\circ \cdots \circ [p]$ applied $i$ times. 

For an arbitrary ideal $I$, we denote 
\[
I_{(p)}:=\sum_{i\geq 0} \; \Span(I^{[p]^i}).
\]
One can show that $I_{(p)}$ is a~$p$-ideal of $(\fa,[p])$ (see, e.g., \cite[Prop. 1.3]{SF}). By definition, $I_{(p)}$ is the smallest $p$-ideal containing the ideal $I$. In particular, $\mathfrak{z}(\fa)_{(p)}=\mathfrak{z}(\fa)$, where $\mathfrak{z}(\fa)$ is the center of $\fa$, a~consequence  of the first condition of (\ref{RRR}).

If $I$ is a~$p$-ideal, then the quotient Lie algebra $\fa/I$ has a~$p$-structure defined by 
\[
(a+I)^{[p]}:= a^{[p]}+I \quad \text{ for any $a\in \fa$},
\]
and the natural map $\pi:\fa\rightarrow \fa/I$ is a~$p$-homomorphism.

For every restricted Lie algebra $\fa$, one can construct its \textit{$p$-enveloping algebra} 
\[
u(\fa):=U(\fa)/I,
\]
where $I$ is the ideal generated by the central elements $a^{[p]}-a^p \in U(\fa)$. 



An $\fa$-module $M$ is called \textit{restricted} if
\[
\underbrace{a\cdot\ldots a}_{p\text{~~times}}\cdot m =a^{[p]}\cdot m \quad \text{for all $a\in \fa$ and any $m\in M$.}
\]

A derivation $\mathscr{D}\in \fder(\fa)$ is called \textit{restricted} if 
\begin{equation}
\label{dvres}
\mathscr{D}(a^{[p]})=(\ad_a)^{p-1} (\mathscr{D}(a)) \quad \text{ for all $a\in \fa$.}
\end{equation}
Denote the space of restricted derivations by $\fder^p(\fa)$. Every inner derivation $\ad_a$, where $a \in \fa$, is a~restricted derivation. Set $\mathfrak{out}^p(\fa):=\mathfrak{der}^p(\fa)/\ad_{\fa}$.  

\subsubsection{Restricted Lie algebra cohomology} 
We denote by $\text{H}^n(\fa;M)$ the usual Chevalley-Eilenberg cohomology of the Lie algebra $\fa$ with coefficient in the $\fa$-module $M.$ Following Hochschild \cite{Ho}, the \textit{restricted} cohomology of a~restricted Lie algebra $\fa$ with coefficients in a~restricted module $M$ is given by
\[
\text{H}^n_{\text{res}}(\fa;M):=\text{Ext}^n_{u(\fa)} (\mathbb{K}, M), \quad \text{where $n\geq0$.} 
\]
Hochschild (\cite{Ho}) showed that there is an exact sequence 
\begin{equation}
\label{six-term}
0\rightarrow \text{H}^1_{\text{res}}(\fa;M) \rightarrow \text{H}^1(\fa;M) \rightarrow S(\fa;M^\fa)\rightarrow \text{H}^2_{\text{res}}(\fa;M)
\rightarrow \text{H}^2(\fa;M)\rightarrow S(\fa; \text{H}^1(\fa;M)),
\end{equation}
where $S(X,Y)$ is the space of $p$-semi-linear\footnote{A map $f:X\rightarrow Y$ is called $p$-\textit{semi-linear} if $f(x+\lambda y)=f(x)+\lambda^p f(y)$ for all $x,y\in X$ and $\lambda \in \mathbb{K}$.} maps $X \rightarrow Y$, and ${M^\fa:=\{m\in M\; |\; \fa\cdot m=0 \}}$ is the space of $\fa$-invariants.

An explicit description of the space of cochains $C^k(\fa;M)$ for $k\leq 3$ was carried out in \cite{EF}. This description was used to classify extensions of restricted modules and infinitesimal deformations of restricted Lie algebras.

The canonical homomorphism 
\[
\text{H}_{\text{res}}(\fa;M)\rightarrow \text{H}(\fa;M)
\]
isomorphically maps $\text{H}^1_{\text{res}}(\fa;M)$ onto the subspace of $\text{H}^1(\fa;M)$ whose elements are represented by the 1-cocycles which satisfy the relation 
\[
x^{p-1}\cdot f(x)=f(x^{[p]_\fa}),
\]
see \cite[Theorem 2.1, page 563]{Ho}. In particular, $\text{H}^1_{\text{res}}(\fa;\fa)\simeq \text{H}^1(\fa;\fa)$ if $\mathfrak{z}(\fa)=0$.

\subsubsection{Restricted outer derivations}\label{Rod}
Clearly, see \cite{Ho, EF}, 
\[
\mathfrak{out}^p(\fa)\simeq \text{H}^1_{\text{res}}(\fa;\fa).
\]
In all the examples we provide in \S \ref{sec: exa}, we do have $\mathfrak{out}^p(\fa)\not =0$. Recall that $\mathfrak{out}(\fg)\simeq\text{H}^1(\fa;\fa) $ while $\mathfrak{out}^p(\fg)\simeq \text{H}^1_{\text{res}}(\fa;\fa)$. For the simple Lie algebras, 
we have $\text{H}^1(\fa;\fa) = \text{H}^1_{\text{res}}(\fa;\fa)$. For Lie algebras with center, such as the Manin triple in \S \ref{Maninp=2}, we have to compute $\text{H}^1_{\text{res}}(\fa;\fa)$ in order to capture restricted outer derivations.

We say that $\mathscr{D}\in \fder^p(\fa)$ has \textit{$p$-property} if there exist $\gamma\in \mathbb{K}$ and $a_0\in \fa$ such that   
\begin{equation}
\label{Conp}
\begin{array}{ccl}
\mathscr{D}^p&=&\gamma \mathscr{D}+\ad_{a_0},\text{ (or, equivalently, $\mathscr{D}^p\simeq \gamma \mathscr{D}$ in $\mathrm{H}^1_{\text{res}}(\fa;\fa)$)}\\
\mathscr{D}(a_0)&=&0.
\end{array}
\end{equation}
  
The restricted outer derivations we provide as main examples in \S \ref{sec: exa} do have $p$-property, same as  nilpotent Lie algebras, see \cite{FSW}:
\begin{equation}\label{FSW}
\begin{array}{l}

\text{Any restricted derivation of any torus identically vanishes (\cite[Prop. 3.1]{FSW})};\\[2pt]

\text{Every restricted derivation of $\mathfrak{hei}(2)$ for $p=2$ is inner\footnote{Recall that the Heisenberg Lie algebra $\mathfrak{hei}(2n)$ is spanned by elements $c_i, a_i$ for $i=1,\ldots,n$, and a~central element $z$ with the only non-zero relations $[c_i,a_i]=z$.} (\cite[Prop. 3.2]{FSW})};\\[2pt]

\text{Apart from a~torus, and $\mathfrak{hei}(2)$ for $p=2$, every outer restricted derivation $\mathscr{D}$}\\ 
\text{of any \textit{nilpotent} restricted Lie algebra satisfies $\mathscr{D}^2=0$ (\cite[Theorem 3.3]{FSW}).}\\

\end{array}
\end{equation}

\subsection{Restricted Lie superalgebras}
%
%
Let $\fa$ be a~ Lie superalgebra defined over a~field of characteristic $p>2$. We denote the parity of a given non-zero homogenous element $a\in\fa$ by $p(a)$; no  confusion with the characteristic of the ground field is possible.

Following \cite{P}, we say that $\fa$ has a~\textit{$p|2p$-structure} if $\fa_\ev$ is restricted and 
\[
\text{$\ad_{a^{[p]}}(b)=(\ad_a)^p(b)$  for all $a \in \fa_\ev$ and $b\in \fa$.}
\]

Recall that the bracket of two odd elements of the Lie superalgebra is polarization of squaring $a\mapsto a^2$. We set 
\[
[2p]:\fa_\od \rightarrow \fa_\ev, \quad a\mapsto (a^2)^{[p]}  \text{~~for any $a\in \fa_\od$}. 
\]
The pair $(\fa, [p|2p])$ is referred to as a~\textit{restricted} Lie superalgebra. 



The following theorem is a~straightforward superization of Jacobson's theorem \ref{Jac}.
\sssbegin{Theorem}\label{SJac}
Let $(e_j)_{j\in J}$ be a~basis of $\fa_\ev$, let the elements $f_j\in \fa_\ev$ be such that  ${(\ad_{e_j})^p=\ad_{f_j}}$. Then, there exists exactly one $p|2p$-mapping $[p|2p]:\fa\rightarrow \fa$ such that 
\[
e_j^{[p]}=f_j \quad \text{ for all $j\in J$}.
\]
\end{Theorem}
\subsubsection{Restricted derivations}\label{RodS}
A derivation $\mathscr{D}\in \fder(\fa)$ is called \textit{restricted} if 
\[
\mathscr{D}(a^{[p]})=(\ad_a)^{p-1} (\mathscr{D}(a))  \quad \text{ for all $a\in \fa_\ev$.}
\]
Consequently, we have 
\[
\mathscr{D}(a^{[2p]})=(\ad_{a^2})^{p-1} (\mathscr{D}(a^2))  \quad \text{ for all $a\in \fa_\od$.}
\]
As in the non-supper setting, the space of restricted derivation is denoted by $\fder^p(\fa)$. 

A~super version of condition (\ref{Conp}):  $\mathscr{D}\in \fder^p_\ev(\fa)$ has \textit{$p$-property} if there exist  $\gamma \in \mathbb{K}$ and $a_0\in \fa_\ev$ such that 
\begin{equation}
\label{SConp}
\mathscr{D}^p=\gamma \mathscr{D}+\ad_{a_0},\ \
\mathscr{D}(a_0)=0.
\end{equation}

\section{The main results}\label{MainR}

\textbf{Let  $\mathscr{K}$ and  $\mathscr{K}^*$ be  $1$-dimensional vector spaces spanned by the vectors $x$ and $x^*$, respectively}.
 
\subsection{The case of Lie algebras}  
The following theorem was proved in \cite{MR} for $\mathbb{K}=\mathbb{R}$. Passing to a~field of characteristic $p\not =2$, the proof is absolutely the same. The case where $p=2$ has been completely described in \cite{BeBou}.

\sssbegin{Theorem}[\cite{MR, BeBou}]\label{MainExt}  Let $\mathscr{B}_\fa$ be $\mathscr{D}$-invariant, where ${\mathscr{D}\in \fder(\fa)}$. 
Then, there exists a~NIS-algebra structure on $\fg:=\mathscr{K} \oplus \fa \oplus \mathscr{K}^*$, where the bracket is defined by
\begin{equation}
\begin{array}{lll}
[\mathscr{K} ,\fg]_\fg&:=&0, \qquad [a,b]_\fg:=[a,b]_\fa+\mathscr{B}_\fa(\mathscr{\mathscr{D}}(a),b)x \text{ for any } a,b\in  \fa,\\[2mm]
[x^*,a]_\fg&:=&\mathscr{D}(a)\text{ for any } a\in  \fa.
\end{array}
\end{equation}
The NIS form $\mathscr{B}_\fg$ on $\fg$ is defined as follows: 
\begin{equation}\label{Leit}
\begin{array}{l}
{ \mathscr{B}_\fg}\vert_{\fa \times \fa}:= \mathscr{B}_\fa, \quad  \mathscr{B}_\fg(x,x^*):=1, \quad  \mathscr{B}_\fg(\fa,x )=\mathscr{B}_\fg(\fa,x^*)=\mathscr{B}_\fg(x,x):=0,  \\ [2mm]\mathscr{B}_\fg(x^*,x^*):=\left\{
\begin{array}{lcl}
\text{arbitrary}, & \text{ if } & p=2\\[2mm]
0,& \text{ if } & p>2\\
\end{array}
\right.
\end{array}
\end{equation}
\end{Theorem}

We call the Lie algebra $(\fg,  \mathscr{B}_\fg)$ constructed in Theorem \ref{MainExt} a~\textit{$\mathscr{D}$-extension} of $(\fa,  \mathscr{B}_{\fa})$  by means of $\mathscr{D}$.

\parbegin{Remark} {\rm If the derivation $\mathscr{D}$ is inner, the double extension is isomorphic to $\fa\oplus\fz$, where $\fz$ is a~$2$-dimensional center, see Theorem $\ref{Isom1}$.}

\end{Remark}

The converse of Theorem \ref{MainExt} is given by the following.

\sssbegin{Theorem}[\cite{MR}] \label{RecMR} If $\mathfrak{z}(\fg)\not =0$, then $(\fg,\mathscr{B}_\fg)$ can be obtained from a~NIS-algebra by means of a~$\mathscr{D}$-extension.
\end{Theorem}

Denote by $\sigma_i^\fa(a,b)$ the coefficient obtained from the expansion
\[
\mathscr{B}_\fa(\mathscr{D}(\lambda a+b),(\ad_{\lambda a+b}^\fa)^{p-2}(a))=\sum_{1\leq i \leq p-1} i\sigma_i ^{\fa}(a,b) \lambda^{i-1}.
\]
For instance, 
\begin{itemize}
\item If $p=2$, then $\sigma_1^{\fa}(a,b)=\mathscr{B}_\fa(\mathscr{D}(b),a)$.
\item If $p=3$, then $\sigma_1^{\fa}(a,b)=\mathscr{B}_\fa(\mathscr{D}(b),[b,a])$ and $\sigma_2^{\fa}(a,b)=2\mathscr{B}_\fa(\mathscr{D}(a),[b,a])$.
\end{itemize}

\parbegin{Lemma}\label{Lem1}
Under assumptions of Theorem \textup{\ref{MainExt}}, we have 
\[
s_i^\fg(a,b)=s_i^\fa(a,b)+ \sigma_i^{\fa}(a,b)x \quad \text{for all $a,b\in \fa$}.
\]
\end{Lemma}

\begin{proof} 
Indeed, 
\[
(\ad_{\lambda a+b}^\fg)^{p-1}(a)=(\ad_{\lambda a+b}^\fa)^{p-1}(a)+\mathscr{B}_\fa(\mathscr{D}(\lambda a+b),(\ad_{\lambda a+b}^\fa)^{p-2}(a) )x.\qed
\]  
\phantom\qedhere 
 \end{proof} 
\subsubsection{The $p$-structures on double extensions. The map $\mathscr{P}$}
Note that for $p=2$ the map $\mathscr{P}$ is quadratic, and therefore it was denoted by $\fq$ in \cite{BeBou, BLS}.

\parbegin{Theorem}[$p$-structure on DE]\label{MainTh} Let $\mathscr{B}_\fa$ be $\mathscr{D}$-invariant, where ${\mathscr{D}\in \fder^p(\fa)}$ and let $\mathscr{D}$ have the $p$-property. For arbitrary $m, l \in \mathbb{K} $, the $p$-structure on $\fa$ can be extended to its $\mathscr{D}$-extension $\fg=\mathscr{K} \oplus \fa \oplus \mathscr{K}^*$ as follows \textup{(}for any $a\in \fa$\textup{)}: 
\[ 
\begin{array}{lcl}
a^{[p]_\fg}&:=&a^{[p]_\fa}+\mathscr{P}(a) x,  \\[2mm]
(x^*)^{[p]_\fg}&:=& a_0+ l x + \gamma x^*,\\[2mm]
x^{[p]_\fg}&:=&m x+b_0,
\end{array}
\]
where $a_0$ and $\gamma$ are as in \eqref{Conp} (the $p$-property), $b_0\in \fz(\fa)$ such that $\mathscr{D}(b_0)=0$, and $\mathscr{P}$ is a~map satisfying (for any $a,b\in \fa$ and any $\lambda\in \mathbb{K} $):
\begin{equation}\label{Qua0}
\begin{array}{rcl}
\mathscr{P}(\lambda a)&=& \lambda^p \mathscr{P}(a),\\[2mm]
\mathscr{P}(a+b)-\mathscr{P}(a)-\mathscr{P}(b)&=&\displaystyle \sum_{1\leq i\leq p-1} \sigma_i^\fa(a,b).
\end{array}
\end{equation}
\end{Theorem}
\begin{proof}
%
Using Jacobson's Theorem \ref{Jac}, it suffices to show that 
\[
\begin{array}{c}
(\ad_a^\fg)^p=\ad_{a^{[p]_\fa}+\mathscr{P}(a) x}, \quad 
(\ad_{x^*}^\fg)^p=\ad_{a_0+ l x +\gamma x^*}, \quad 
(\ad_{x}^\fg)^p=\ad_{m x+b_0}.
\end{array}
\]
Indeed, let $f=ux+b+vx^*\in \fg$. We have

\[
\begin{array}{lcl}
\ad_{mx+b_0}^\fg(f)-(\ad_{x}^\fg)^p (f)&=&[mx+b_0,f]_\fg-(\ad_{x}^\fg)^p (f) \\[2mm]
&=&[b_0,b]_\fa+\mathscr{B}_\fa(\mathscr{D}(b_0),b)x-v\mathscr{D}(b_0)=0,
\end{array}
\]
since $x$ and $b_0$ are central, and $\mathscr{D}(b_0)=0$. Besides, using condition (\ref{Conp}) we get 

\[
\begin{array}{lcl}
\ad_{a_0+l x + \gamma x^*}^\fg(f)-(\ad_{x^*}^\fg)^p (f)&=&[ a_0+l x + \gamma x^*,f]_\fg-\mathscr{D}^p (b) \\[2mm]
&=&[a_0,b]_\fa+\mathscr{B}_\fa(\mathscr{D}(a_0),b)x-v\mathscr{D}(a_0)+\gamma \mathscr{D}(b)-\mathscr{D}^p (b)\\[2mm]
&=&0,
\end{array}
\]
and
\[
\begin{array}{lcl}
\ad_{a^{[p]_\fg}  }^\fg(f)-(\ad_{a}^\fg)^p (f)&=&[ a^{[p]_\fa}+\mathscr{P}(a) x ,f]_\fg-(\ad_{a}^\fg)^p (b)+v(\ad_{a}^\fg)^{p-1} (\mathscr{D}(a)) \\[2mm]
&=&[ a^{[p]_\fa}, b]_\fa+\mathscr{B}_\fa(\mathscr{D}(a^{[p]_\fa}),b)x -v \mathscr{D}(a^{[p]_\fa})-(\ad_{a}^\fg)^p (b)\\[2mm]
&&+v(\ad_{a}^\fg)^{p-1} (\mathscr{D}(a)) \\[2mm]
&=&[ a^{[p]_\fa}, b]_\fa+\mathscr{B}_\fa(\mathscr{D}(a^{[p]_\fa}),b)x -v \mathscr{D}(a^{[p]_\fa})-(\ad_{a}^\fa)^p (b)\\[2mm]
&&-\mathscr{B}_\fa(\mathscr{D}(a),(\ad_a)^{p-1}(b))x+v(\ad_{a}^\fa)^{p-1} (\mathscr{D}(a))\\[2mm]
&&+v\mathscr{B}_\fa(\mathscr{D}(a),(\ad_{a}^\fg)^{p-2} (\mathscr{D}(a)))x \\[2mm]
&=&0,

\end{array}
\]
because
\[
\mathscr{B}_\fa(\mathscr{D}(a),(\ad_a)^{p-1}(b))=(-1)^{p-1}\mathscr{B}_\fa((\ad_a)^{p-1}\circ \mathscr{D}(a),b),
\]
and also (using again the $\mathscr{D}$-invariance and the $p$-property):
\[
\begin{array}{lcl}
\mathscr{B}_\fa(\mathscr{D}(a),(\ad_a^\fa)^{p-2}(\mathscr{D}(a)))&=&
\left \{
\begin{array}{ll}
\mathscr{B}_\fa(\mathscr{D}(a), (\ad_a^\fa)^{\frac{p-3}{2}} \circ \ad_a^\fa\circ (\ad_a^\fa)^{\frac{p-3}{2}} (\mathscr{D}(a)))&\text{if $p>2$}\\[2mm]
\mathscr{B}_\fa(\mathscr{D}(a),\mathscr{D}(a))& \text{if $p=2$}
\end{array}
\right.
\\[3mm]
&=&
\left \{
\begin{array}{ll}
(-1)^{\frac{p-3}{2}}\mathscr{B}_\fa((\ad_a)^{\frac{p-3}{2}} (\mathscr{D}(a)),\ad_a \circ (\ad_a)^{\frac{p-3}{2}} (\mathscr{D}(a)))& \text{if $p>2$}\\[2mm]
\mathscr{B}_\fa(\gamma \mathscr{D}(a)+[a_0,a],a)& \text{if $p=2$}
\end{array}
\right .
\\[3mm]
&=&0.
\end{array}
\]
Now, let $a,b \in \fa$. Using Lemma \ref{Lem1}, we have
\[
\begin{array}{lcl}
(a+b)^{[p]_\fg}&=&(a+b)^{[p]_\fa} +\mathscr{P}(a+b)x\\[2mm]
&=& a^{[p]_\fa} + b^{[p]_\fa}+ \displaystyle \sum_{1\leq i\leq p-1}s_i^\fa(a,b) +\mathscr{P}(a+b)x\\[2mm]
&=& a^{[p]_\fg} + b^{[p]_\fg}+ \displaystyle \sum_{1\leq i\leq p-1}s_i^\fg(a,b) +\mathscr{P}(a+b)x-\mathscr{P}(a)x-\mathscr{P}(b)x-\sum_{1\leq i\leq p-1}\sigma_i^{\fa}(a,b) x.\\[2mm]
\end{array}
\]
By virtue of Eq. (\ref{Qua0}), it follows that $[p]_\fg$ is a~$p$-mapping.
 
The proof is now complete.
\end{proof}

The converse of Theorem \ref{MainTh} is the following.

\parbegin{Theorem}[Converse of Theorem \ref{MainTh}]
\label{Rec1}
Let $(\fg,\mathscr{B}_\fg)$ be an irreducible and restricted NIS-algebra of $\text{dim}>1$ such that $\fz(\fg)\not=0$. Let $0\not=x\in \fz(\fg)$ be such that $\mathscr{K}^\perp$ is a~$p$-ideal. Then, $(\fg,\mathscr{B}_\fg)$ is a~$\mathscr{D}$-extension of a~restricted NIS-algebra $(\fa,\mathscr{B}_\fa)$, where $\mathscr{D}$ is a~restricted derivation with $p$-property.
\end{Theorem}
\begin{proof}
The subspace $\mathscr{K}$ is an ideal in $(\fg,\mathscr{B}_\fg)$ because $x$ is central in $\fg$. Moreover, $\mathscr{K}^\perp$ is also an ideal in $(\fg,\mathscr{B}_\fg)$, see \cite{BeBou, MR}.  Since $\fg$ is irreducible, it follows that $\mathscr{K}\subset \mathscr{K}^\perp$ and $\dim(\mathscr{K}^\perp)=\dim(\fg)-1$. Therefore, there exists a~non-zero $x^* \in \fg$ such that
\[
\fg=\mathscr{K}^\perp\oplus \mathscr{K}^*. 
\]
We normalize $x^*$ so that $\mathscr{B}_\fg(x,x^*)=1$. Besides, $\mathscr{B}_\fg(x,x)=0$ since ${\mathscr{K}\cap \mathscr{K}^\perp=\mathscr{K}}$.

Set $\fa:=(\mathscr{K} +\mathscr{K}^*)^\perp$. We then obtain a~decomposition $\fg=\mathscr{K} \oplus \fa \oplus \mathscr{K}^*$. 

There exists a~NIS-algebra structure on the vector space $\fa$ for which $\fg$ is its double extension by Theorem \ref{RecMR}. We denote a~NIS on $\fa$ by $\mathscr{B}_\fa$.

It remains to show that there is a~$p$-structure on $\fa$. Since $\fa \subset \mathscr{K}^\perp$, then 
\[
a^{[p]_\fg}\in \mathscr{K}_p^\perp=\mathscr{K}^\perp=\mathscr{K} \oplus \fa \quad \text{ for any $a \in \fa$}.
\]  
It follows that
\begin{equation*}
a^{[p]_\fg}=\mathscr{P}(a)x+s(a), \quad \text{where $s(a)\in \fa$}.
\end{equation*}

Define a~$p$-structure on $\fa$ by setting 
\[
s:\fa \rightarrow \fa, \quad a\mapsto s(a).
\]

The fact that $(\lambda a)^{[p]_\fg}=\lambda^p(a^{[p]_\fg})$ implies $s(\lambda a)=\lambda^p s(a)$ and $\mathscr{P}(\lambda a)=\lambda^p \mathscr{P}(a)$. Besides,
\[
\begin{array}{lcl}
0&=&[a^{[p]_\fg},b]_\fg-(\ad_a^\fg)^{p}(b)\\[2mm]
&=&[s(a),b]_\fa+\mathscr{B}_\fa(\mathscr{D}(s(a)), b)x-(\ad_a^\fa)^p(b)-\mathscr{B}_\fa(\mathscr{D}(a),(\ad_a^\fa)^{p-1}(b))x.
\end{array}
\]
Therefore, 
\[
\mathscr{D}(s(a))= (\ad_a^\fa)^{p-1}\circ (\mathscr{\mathscr{D}}(a)), \quad [s(a),b]_\fa=(\ad_a^\fa)^p(b).
\]
Now, since
\[
\begin{array}{lcl}
\displaystyle \sum_{1\leq i \leq p-1} is_i^\fg(a,b) \lambda^{i-1}&=&(\ad_{\lambda a+b}^\fg)^{p-1}(a),\\[2mm]
&=&(\ad_{\lambda a+b}^\fa)^{p-1}(a)+\mathscr{B}_\fa(\mathscr{\mathscr{D}}(\lambda a+b),(\ad_{\lambda a+b}^\fa)^{p-2}(a) )x,
\end{array}
\]
it follows that 
\[
s^\fg_i(a,b)=s_i^{\fa}(a,b)+\sigma_i^{\fa}(a,b)x.
\]
Moreover, 
\[
\begin{array}{lcl}
0&=&(a+b)^{[p]_\fg}-a^{[p]_\fg}-b^{[p]_\fg}-\displaystyle \sum_{1\leq i \leq p-1}s_i^\fg(a,b)\\[2mm]
&=&(\mathscr{P}(a+b)-\mathscr{P}(a)-\mathscr{P}(b)- \displaystyle \sum_{1\leq i \leq p-1}\sigma_i^{\fa}(a,b) ) x-\displaystyle \sum_{1\leq i \leq p-1}s_i^{\fa}(a,b)\\
&&+s(a+b)-s(a)-s(b).
\end{array}
\]
Consequently, 
\begin{eqnarray}
\nonumber s(a+b)&=&s(a)+s(b)+\displaystyle\sum_{1\leq i \leq p-1}s_i^{\fa}(a,b), \\[2mm]
\nonumber \label{condq2} \mathscr{P}(a+b)-\mathscr{P}(a)-\mathscr{P}(b)&=&\displaystyle \sum_{1\leq i \leq p-1}\sigma_i^{\fa}(a,b).
 \end{eqnarray}
It follows that $s$ defines a~$p$-mapping on $\fa$, that $\mathscr{D}$ is a~restricted derivation of $\fa$ (relative to the $p$-mapping $s$), and $\mathscr{P}$ is a~mapping on $\fa$ satisfying Eq. (\ref{Qua0}).

Suppose that 
\begin{equation}
\label{xsp}
(x^*)^{[p]_\fg}=a_0+\beta x+\gamma x^*, \text{ where $a_0\in \fa$ and $\beta, \gamma \in \mathbb{K} $}.
\end{equation}
For all $a\in \fa$, we have
\[
0=[(x^*)^{[p]},a]_\fg-(\ad_{x^*}^\fg)^p(a)=[a_0,a]+ \mathscr{B}_\fa(\mathscr{D}(a_0),a)x+\gamma \mathscr{D}(a)-\mathscr{D}^{p}(a).
\]
Since $\mathscr{B}_\fg$ is non-degenerate, it follows that $\mathscr{D}^p=\gamma \mathscr{D} +\ad_{a_0}$ and $\mathscr{D}(a_0)=0$.

Suppose now that 
\begin{equation}
\label{xp}
x^{[p]_\fg}=b_0+mx+\delta x^* \text{ for some $m,\delta\in \mathbb{K} $ and $b_0\in \fa$}.
\end{equation}

For any $b\in \fa$, we have
\[
0=[x^{[p]_\fg},b]_\fg-(\ad_x^\fg)^{p}(b)=[b_0,b]_\fa+\mathscr{B}_\fa(\mathscr{D}(b_0), b)x + \delta \mathscr{D}(b). 
\]
Since $\mathscr{B}_\fa$ is non-degenerate, it follows that $\delta \mathscr{D}(b)+[b_0,b]_\fa=0$ and $\mathscr{D}(b_0)=0$.

\underline{The case where $\mathscr{D}\not\simeq0$ in $\mathrm{H}^1_{\mathrm{res}}(\fa;\fa)$:} It follows that $\delta=0$ and $b_0\in \fz(\fa)$. Therefore, $\fg$ can be obtained from the restricted Lie algebra $\fa$ by means of the restricted derivation $\mathscr{D}$ as in Theorem \ref{MainTh}

\underline{The case where $\mathscr{D}\simeq0$ in $\mathrm{H}^1_{\mathrm{res}}(\fa;\fa)$:} Without loss of generality, we can assume that ${\mathscr{D}=0}$, cf. Theorem \ref{Isom1}. In this case, $\fg \simeq\fa\oplus \fz$, where $\fz$ is the center of $\fg$ spanned by $x$ and $x^*$, and the $p$-structure is given by Eqs. (\ref{xsp}), (\ref{xp}), where $a_0,b_0\in  \fz(\fa)$, while $\gamma, \delta, m$ and $\beta$ are arbitrary. \end{proof}


\subsection{The case of Lie superalgebras: $p(\mathscr{B}_\fa)=\ev$}\label{sec: even}
Let $(\fa, \mathscr{B}_\fa)$ be a~NIS-superalgebra, $p\not=2$. A double extension by means of an even derivation is called a~\textit{$\mathscr{D}_\ev$-extension}, while a~double extension by means of an odd derivation is called a~\textit{$\mathscr{D}_\od$-extension}.

For every $a,b\in \fa_\ev$, denote by $\sigma_i^\fa(a,b)$ the coefficient obtained from the expansion
\[
\mathscr{B}_\fa(\mathscr{D}(\lambda a+b),(\ad_{\lambda a+b}^\fa)^{p-2}(a))=\sum_{1\leq i \leq p-1} i\sigma_i ^{\fa}(a,b) \lambda^{i-1}.
\]

Let $\mathscr{B}_\fa$ be $\mathscr{D}$-invariant, where ${\mathscr{D}\in \fder^p(\fa)}$.  If $p(\mathscr{D} )=\od$, assume that $\mathscr{D} $ satisfies the following 
 conditions for some $b_0\in \fa_\ev$:
\begin{equation}
\label{DB0}
\mathscr{D}(b_0)=0,\ \ 
\mathscr{D}^2=\ad_{b_0}, \ \ 
\mathscr{B}_\fa(b_0,b_0)=0.
\end{equation}
For $p=3$, let us suppose further that\footnote{This condition is automatically satisfied for $p>3$, since $\mathscr D$ is a~derivation.}
\begin{equation}
\label{p=3Con1}
\mathscr{B}_\fa ({\mathscr D}(a), [a,a]_\fa)=0 \quad \text{for all $a\in \fa_\od$}.
\end{equation}

Let  
$p(x)=p(x^*)=p(\mathscr{D})$.  Let $b_0$ be as in Eq. (\ref{DB0}). 
Following \cite{BB, ABB}, there exists a~NIS-superalgebra structure on $\fg:=\mathscr{K} \oplus \fa \oplus \mathscr{K} ^*$, where the  bracket is defined by 
\begin{equation}\label{*}
\begin{array}{lcl}
[\mathscr{K} ,\fg]_\fg&:=&0,\\[2mm]
 [a,b]_\fg & := & [a,b]_\fa+(-1)^{p(\mathscr{D})}\mathscr{B}_\fa(\mathscr{D}(a),b)x \text{ for any } a,b\in  \fa, \\[2mm]
 [x^*,a]_\fg&:=&\mathscr{D}(a) - \left \{ 
  \begin{array}{ll}
 2 \mathscr{B}_\fa(a,b_0)x & \text{if $p(\mathscr{D})=\od$}\\[2mm]
  0 & \text{if $p(\mathscr{D})=\ev$}
  \end{array}
  \text{ for any } a\in  \fa
  \right .
  \\[2mm]
  [x^*,x^*]_\fg&:=& \left \{ 
  \begin{array}{ll}
 2 b_0 & \text{if $p(\mathscr{D})=\od$}\\[2mm]
  0 & \text{if $p(\mathscr{D})=\ev$}
  \end{array}
  \right .
\end{array}
\end{equation}
The NIS form $\mathscr{B}_\fg$ on $\fg$ is defined as in (\ref{Leit}). 

\sssbegin{Theorem}[$p$-structure on DE for $p(\mathscr{B}_\fa)=\ev$]\label{DoBe} Let $\mathscr{B}_\fa$ be $\mathscr{D}$-invariant, where ${\mathscr{D}\in \fder^p(\fa)}$ and let $p(\mathscr{B}_\fa)=\ev$.

\textup{(}i\textup{)} If $p(\mathscr{D})=\od$, $\mathscr{D}$ satisfies the conditions \textup{(\ref{DB0})} and \textup{(\ref{p=3Con1})} and 
\begin{equation}\label{StarT}
2\mathscr{B}_\fa(a^{[p]_\fa}, b_0)-\mathscr{B}_\fa(\mathscr{D}(a), ({\ad_a})^{p-2}\circ \mathscr{D}(a))=0 \text{ for all $a\in \fa_\ev$},
\end{equation}
then the $p|2p$-structure on $\fa$ can be extended to its double extension $\fg$ as follows:
\[ 
\begin{array}{c}
a^{[p]_\fg}:=a^{[p]_\fa} \quad \text{for any $a\in \fa_\ev$}.
\end{array}
\]

\textup{(}ii\textup{)}  If $p(\mathscr{D})=\ev$ and $\mathscr{D}$ has $p$-property, then for arbitrary $m, l \in \mathbb{K}$, the $p|2p$-structure on $\fa$ can be extended to its double extension $\fg$ as follows (for any $a\in \fa_\ev$):
\[ 
\begin{array}{lcl}
a^{[p]_\fg}&:=&a^{[p]_\fa}+\mathscr{P}(a) x,  \\[2mm]
(x^*)^{[p]_\fg}&:=& a_0+ l x + \gamma x^*,\\[2mm]
x^{[p]_\fg}&:=&m x+c_0,
\end{array}
\]
where $\gamma, a_0$ are defined as in \eqref{SConp} (the $p$-property), $c_0\in \fz_\ev(\fa)$ such that $\mathscr{D}(c_0)=0$, and $\mathscr{P}$ is a~map satisfying (for any $a,b\in \fa_\ev$ and for any $\lambda\in \mathbb{K}$):
\begin{equation}\label{Qua}
\begin{array}{rcl}
\mathscr{P}(\lambda a)&=& \lambda^p \mathscr{P}(a),\\[2mm]
\mathscr{P}(a+b)-\mathscr{P}(a)-\mathscr{P}(b)&=&\displaystyle \sum_{1\leq i\leq p-1} \sigma_i^\fa(a,b).
\end{array}
\end{equation}
\end{Theorem}

We need the following two lemmas.
\parbegin{Lemma}\label{SLem1} \textup{(i)} If $p(\mathscr{B}_\fa)=\ev$ and $p(\mathscr{D})=\od$, then
\[
s_i^\fg(a,b)=s_i^\fa(a,b)  \quad \text{for all $a,b\in \fa_\ev$}.
\]
\textup{(ii)} If $p(\mathscr{B}_\fa)=p(\mathscr{D})=\ev$, then
\[
s_i^\fg(a,b)=s_i^\fa(a,b)+ \sigma_i^{\fa}(a,b)x  \quad \text{for all $a,b\in \fa_\ev$}.
\]
\end{Lemma}
\begin{proof}
Indeed, 
\[
(\ad_{\lambda a+b}^\fg)^{p-1}(a)=(\ad_{\lambda a+b}^\fa)^{p-1}(a)+ \left \{
\begin{array}{ll}\mathscr{B}_\fa(\mathscr{D}(\lambda a+b),(\ad_{\lambda a+b}^\fa)^{p-2}(a) )x &\text{ if $p(\mathscr{D})=\ev$}\\[1mm]
0 &\text{ if $p(\mathscr{D})=\od$}
\end{array}
\right..
\]
The result follows immediately. 
\end{proof}

\parbegin{Lemma}\label{SLem2} For all $a\in \fa_\ev$  and for all $b\in \fa$, we have
\[
\mathscr{B}_\fa(\mathscr{D}(a^{[p]_\fa}),b)=\mathscr{B}_\fa(\mathscr{D}(a),(\ad_a)^{p-1}(b))\quad \text{and} \quad  
\mathscr{B}_\fa(\mathscr{D}(a),(\ad_a^\fa)^{p-2}(\mathscr{D}(a)))=0.
\]
\end{Lemma}
\begin{proof}
Indeed, 
\[
\begin{split}
\mathscr{B}_\fa(\mathscr{D}(a),(\ad_a)^{p-1}(b))=(-1)^{p-1}\mathscr{B}_\fa((\ad_a)^{p-1}\circ \mathscr{D}(a),b)=&\mathscr{B}_\fa(\mathscr{D}(a^{[p]_\fa}),b).
\end{split}
\]

Besides, 
\[
\begin{array}{lcl}
\mathscr{B}_\fa(\mathscr{D}(a),(\ad_a^\fa)^{p-2}(\mathscr{D}(a)))&=&
\mathscr{B}_\fa(\mathscr{D}(a), (\ad_a^\fa)^{\frac{p-3}{2}} \circ \ad_a^\fa\circ (\ad_a^\fa)^{\frac{p-3}{2}} (\mathscr{D}(a)))\\[2mm]
&=&
(-1)^{\frac{p-3}{2}}\mathscr{B}_{\fa}((\ad_a)^{\frac{p-3}{2}}\circ \mathscr{D}(a),\ad_a \circ (\ad_a)^{\frac{p-3}{2}}\circ \mathscr{D}(a))=0.
\qed
\end{array}
\]
\phantom\qedhere 
\end{proof}
Let us prove Theorem \ref{DoBe}.
\begin{proof} (i) Using Jacobson's Theorem \ref{SJac} in super setting, it is enough to show that 
\[
\begin{array}{c}
(\ad_a^\fg)^p=\ad_{a^{[p]_\fa}}.
\end{array}
\]
Indeed, let $f=ux+b+vx^*\in \fg$. 
\[
\begin{array}{lcl}
\ad_{a^{[p]_\fa}}^\fg(f)-(\ad_{a}^\fg)^p (f)&=&[ a^{[p]_\fa} ,b]_\fg+v[a^{[p]_\fa},x^*]_\fg-(\ad_{a}^\fg)^p (b)-v(\ad_{a}^\fg)^{p} (x^*) \\[2mm]
&=&[ a^{[p]_\fa}, b]_\fa-\mathscr{B}_\fa(\mathscr{D}(a^{[p]_\fa}),b)x-v \mathscr{D}(a^{[a]_\fa})+2v \mathscr{B}_\fa(a^{[p]_\fa},b_0)x \\[2mm]
&&-(\ad_{a}^\fg)^p (b)+v(\ad_{a}^\fg)^{p-1} (\mathscr{D}(a)-2 \mathscr{B}_\fa(a,b_0)x) \\[2mm]
&=&[ a^{[p]_\fa}, b]_\fa-\mathscr{B}_\fa(\mathscr{D}(a^{[p]_\fa}),b)x -v \mathscr{D}(a^{[p]_\fa})+2 v \mathscr{B}_\fa(a^{[p]_\fa},b_0)x\\[2mm]
&&-(\ad_{a}^\fa)^p (b)+\mathscr{B}_\fa(\mathscr{D}(a),(\ad_a)^{p-1}(b))x\\[2mm]
&&+v(\ad_{a}^\fa)^{p-1} (\mathscr{D}(a))-v\mathscr{B}_\fa(\mathscr{D}(a),(\ad_{a}^\fg)^{p-2} (\mathscr{D}(a)))x. \\[2mm]
&=&0,

\end{array}
\]
because of Lemma \ref{SLem2} and the fact that $2\mathscr{B}_\fa(a^{[p]_\fa},b_0)-\mathscr{B}_\fa(\mathscr{D}(a),(\ad_{a}^\fg)^{p-2} (\mathscr{D}(a)))=0$.

Since $s_i^\fa(a,b)=s_i^\fg(a,b)$, it follows that
\[
(a+b)^{[p]_\fg}-a^{[p]_\fg}-b^{[p]_\fg}-\sum_{1\leq i \leq p-1}s_i^\fg(a,b)=(a+b)^{[p]_\fa}-a^{[p]_\fa}-b^{[p]_\fa}-\sum_{1\leq i \leq p-1}s_i^\fa(a,b)=0.
\]
(ii) Using Jacobson's Theorem \ref{SJac}, it is enough to show that 
\[
\begin{array}{c}
(\ad_a^\fg)^p=\ad_{a^{[p]_\fa}+\mathscr{P}(a) x}, \quad 
(\ad_{x^*}^\fg)^p=\ad_{a_0+ l x +\gamma x^*}, \quad 
(\ad_{x}^\fg)^p=\ad_{m x+c_0}.
\end{array}
\]
Indeed, let $f=ux+b+vx^*\in \fg$. We have

\[
\begin{array}{lcl}
\ad_{mx+c_0}^\fg(f)-(\ad_{x}^\fg)^p (f)&=&[mx+c_0,f]_\fg-(\ad_{x}^\fg)^p (f) \\[2mm]
&=&[c_0,b]_\fa+\mathscr{B}_\fa(\mathscr{D}(c_0),b)x-v\mathscr{D}(c_0)=0,
\end{array}
\]
since $x\in \fz_\ev(\fg)$, $c_0\in \fz_\ev(\fa)$ and $\mathscr{D}(c_0)=0$. Besides, using condition (\ref{Conp}) we obtain  

\[
\begin{array}{lcl}
\ad_{a_0+l x + \gamma x^*}^\fg(f)-(\ad_{x^*}^\fg)^p (f)&=&[ a_0+l x + \gamma x^*,f]_\fg-\mathscr{D}^p (b) \\[2mm]
&=&[a_0,b]_\fa+\mathscr{B}_\fa(\mathscr{D}(a_0),b)x-v\mathscr{D}(a_0)+\gamma \mathscr{\mathscr{D}}(b)-\mathscr{D}^p (b)\\[2mm]
&=&0.
\end{array}
\]
With Lemma \ref{SLem2}, we show by a~direct computation that $\ad_{a^{[p]_\fa}+\mathscr{P}(a) x}^\fg-(\ad_{a}^\fg)^p=0$. 

Now, let $a,b \in \fa_\ev$. Using Lemma \ref{SLem1}, we have
\[
\begin{array}{lcl}
(a+b)^{[p]_\fg}&=&(a+b)^{[p]_\fa} +\mathscr{P}(a+b) x\\[2mm]
&=& a^{[p]_\fa} + b^{[p]_\fa}+ \displaystyle \sum_{1\leq i\leq p-1}s_i^\fa(a,b) +\mathscr{P}(a+b)x\\[2mm]
&=& a^{[p]_\fg} + b^{[p]_\fg}+ \displaystyle \sum_{1\leq i\leq p-1}s_i^\fg(a,b) +  ( \mathscr{P}(a+b)-\mathscr{P}(a)-\mathscr{P}(b) - \displaystyle \sum_{1\leq i\leq p-1}\sigma_i^{\fa}(a,b)  )x.\\[2mm]
\end{array}
\]
By Eq. (\ref{Qua}), it follows that $[p]_\fg$ defines a~$p$-structure on $\fg$.

The proof is now complete.
\end{proof}

\sssbegin{Theorem}[Converse of Theorem \ref{DoBe}]
\label{DeBeRec}
Let $(\fg,\mathscr{B}_\fg)$ be an irreducible and restricted NIS-superalgebra of $\text{dim}>~1$.

\textup{(}i\textup{)} Suppose there exists $0\not=x\in \fz_\ev(\fg)$ such that $\mathscr{K}^\perp$ is a~$p$-ideal. 
Then, $(\fg,\mathscr{B}_\fg)$ is a~$\mathscr{D}_\ev$-extension of a~restricted NIS-superalgebra $(\fa,\mathscr{B}_\fa)$, where $\mathscr{D}$ is a~restricted derivation with $p$-property.

\textup{(}ii\textup{)}  Let $\fz_\od(\fg)\not =0$. Then, $(\fg,\mathscr{B}_\fg)$ is a~$\mathscr{D}_\od$-extension of a~restricted NIS-superalgebra $(\fa,\mathscr{B}_\fa)$, where $\mathscr{D}_\od$ is a restricted derivation that satisfies conditions \textup{(\ref{DB0})} and \textup{(\ref{p=3Con1})}.

\end{Theorem}

\begin{proof} 
Similar to that of Theorem \ref{Rec1}, minding the Sign Rule.
\end{proof}

\subsection{The case of Lie superalgebras: $p(\mathscr{B}_\fa)=\od$}
Let $(\fa, \mathscr{B}_\fa)$ be a~NIS-Lie superalgebra, $p\not =2 $. Let $ \mathscr{B}_\fa$ be $\mathscr{D}$-invariant for some $\mathscr{D}\in \fder^p(\fa)$ that satisfies the following conditions:
\begin{align}
\label{BoD1} &\text{If $p(\mathscr{D})=\od$, we assume that }
\mathscr{D}^2=\ad_{b_0} \text{ and }
 \mathscr{D}(b_0)=0.
\\
&  \label{p=3Con2} \text{If $p(\mathscr{D})=\ev$ and $p=3$, we assume that } 
 \mathscr{B}_\fa (\mathscr{D}(a),[a,a]_\fa)=0\quad \text{for any $a \in \fa_\od$}.
\end{align}

Let ${\fg:=\mathscr{K} \oplus \fa \oplus \mathscr{K}^*}$, where 
$p(x)=p(\mathscr{D})+\od$ and 
$p(x^*)=p(\mathscr{D})$. As shown in \cite{ABBQ}, there exists a~NIS-superalgebra structure on $\fg$. Let $\lambda_0\in \mathbb{K}$ and $b_0$ be as in Eq.~(\ref{BoD1}). Set
\begin{equation}
\label{lambda}
\begin{array}{lcl}
[\mathscr{K} ,\fg]_\fg&:=&0,\\[2mm]
[a,b]_\fg&:=&[a,b]_\fa+\mathscr{B}_\fa(\mathscr{D}(a),b)x \text{ for any $a,b\in \fa$} \\[2mm]
[x^*,x^*]_\fg&:=& 
\left \{
\begin{array}{ll} 
2 b_0+\lambda_0x, & \text{if $p(\mathscr{D})=\od$}   \\[2mm]
0, & \text{if $p(\mathscr{D})=\ev$} 
\end{array}
\right. ,\\[2mm]
 [x^*,a]_\fg&:=&\mathscr{D}(a) - \left \{
\begin{array}{ll} 
(-1)^{p(a)}2\mathscr{B}_\fa(a,b_0)x, & \text{if $p(\mathscr{D})=\od$}   \\[2mm]
0, & \text{if $p(\mathscr{D})=\ev$} 
\end{array}
\right.\text{ for any } a\in  \fa.
\end{array} 
\end{equation}
The NIS form $\mathscr{B}_\fg$ on $\fg$ is defined as in (\ref{Leit}). 


\sssbegin{Theorem}[$p$-structure on DE for $p(\mathscr{B}_\fa)=\od$]\label{DoBo} Let $\mathscr{B}_\fa$ be $\mathscr{D}$-invariant, where ${\mathscr{D}\in \fder^p(\fa)}$, and let $p(\mathscr{B}_\fa)=\od.$

\textup{(}i\textup{)}  If $p(\mathscr{D})=\od$ and condition \textup{(\ref{BoD1})} is satisfied, then the $p|2p$-structure on $\fa$ can be extended to $\fg$ as follows \textup{(for any $a\in \fa_\ev$, $m\in \mathbb K$ and $c_0\in \fz_\ev(\fa)$  such that $\mathscr{D}(c_0)=0$)}
\[ 
\begin{array}{c}
a^{[p]_\fg}:=a^{[p]_\fa}+\mathscr{P}(a) x,  \quad 
x^{[p]_\fg}:=m x+c_0.
\end{array}
\]
where $\mathscr{P}$ is a~mapping on $\fa_{\ev}$ satisfying (for any $a,b\in \fa_\ev$ and for any $\lambda\in \mathbb{K}$):
\begin{equation}\label{Qua2}
\begin{array}{rcl}
\mathscr{P}(\lambda a)&=& \lambda^p \mathscr{P}(a),\\[2mm]
\mathscr{P}(a+b)-\mathscr{P}(a)-\mathscr{P}(b)&=&\displaystyle \sum_{1\leq i\leq p-1} \sigma_i^\fa(a,b).
\end{array}
\end{equation}

\textup{(}ii\textup{)}  If $p(\mathscr{D})=\ev$, has $p$-property, and satisfies  condition \textup{(\ref{p=3Con2})}, then the $p|2p$-structure on $\fa$ can be extended to $\fg$ as follows \textup{(for any $a\in \fa_\ev$)}
\[ 
\begin{array}{c}
a^{[p]_\fg}:=a^{[p]_\fa},  \quad 
(x^*)^{[p]_\fg}:= \gamma x^*+a_0,
\end{array}
\]
where $a_0$ and $\gamma$ are as in the $p$-property. 
\end{Theorem}

\begin{proof}
Similar to that of Theorem \ref{DoBe}.
\end{proof}

\sssbegin{Theorem}[Converse of Theorem \ref{DoBo}]
\label{DoBoRec}
Let $(\fg,\mathscr{B}_\fg)$ be an irreducible restricted NIS-superalgebra of $\text{dim}>1$ such that $p(\mathscr{B}_\fa)=\od$. 

\textup{(}i\textup{)}  Suppose that $\fz_\ev(\fg) \not =0$. Then, $(\fg,\mathscr{B}_\fg)$ is a~$\mathscr{D}_\od$-extension of a~restricted NIS-superalgebra $(\fa,\mathscr{B}_\fa)$, where $\mathscr{D}_\od$ is an odd restricted derivation.

\textup{(}ii\textup{)}  Suppose there exists $0\not=x\in \fz_\od(\fg)$ such that $\mathscr{K}^\perp$  is a~$p$-ideal. 
Then, $(\fg,\mathscr{B}_\fg)$ is a~$\mathscr{D}_\ev$-extension of a~restricted NIS-superalgebra $(\fa,\mathscr{B}_\fa)$ such that $\mathscr{B}_\fa$ is $\mathscr{D}$-invariant, where $\mathscr{D}_\ev$ is a~restricted derivation with $p$-property and (only for $p=3$) satisfying condition \textup{(\ref{p=3Con2})}.
\end{Theorem}

\begin{proof}
Similar to that of Theorem \ref{Rec1}, minding the Sign Rule.
\end{proof}

\subsection{Isomorphisms of NIS-superalgebras}
For a NIS-superalgebra $\fa$ with NIS $\mathscr{B}_\fa$, denote by $\fg$ (resp. $\tilde \fg$) the double extension of $\fa$ by means of a~derivation $\mathscr{D}$ (resp. $\tilde{\mathscr{D}}$).  An isomorphism $\pi:\fg\rightarrow \tilde \fg$ should satisfy (see \cite{FS,BeBou}): 
\[
\begin{array}{lcll}
\pi([f,g]_\fg)&=&[\pi(f), \pi(g)]_{\tilde \fg} & \text{for any $f, g\in \fg$},\\[2mm]
\mathscr{B}_{\tilde \fg}(\pi(f), \pi(g))&=&\mathscr{B}_\fg(f,g) & \text{for any $f,g\in \fg$}.
\end{array}
\]
We will further assume that $\pi(\mathscr{K} \oplus \fa)=\tilde{\mathscr{K}} \oplus \tilde{\fa}$, and call $\pi$ an \textit{adapted isomorphism}. We will see how the derivations $\mathscr{D}$ and $\tilde {\mathscr{D}}$ are related with each other when $\fg$ and $\tilde \fg$ are isomorphic.

\subsubsection{The case of Lie algebras}
The following theorem was proved in \cite{FS} in the case where $\mathbb{K} =\mathbb{R}$. The passage to any $p\not =2$, the proof is absolutely the same. The case where $p=2$ was studied in \cite{BeBou}.

\parbegin{Theorem}[\cite{FS, BeBou}] \label{Isom1}
Let $\mathscr{B}_\fa$ be $\mathscr{D}$- and $\tilde {\mathscr{D}}$-invariant, where $\mathscr{D}, \tilde {\mathscr{D}} \in \fder(\fa)$. There exists an adapted isomorphism $\pi:  \fg\rightarrow  \tilde \fg$ if and only if there exists an automorphism $\pi_0:\fa\rightarrow \fa,$ a~scalar $\lambda \in \mathbb{K} ^{\times}$ and $\varkappa \in \fa$ \textup{(}all unique\textup{)} such that
\begin{eqnarray}
\label{Cb} \pi^{-1}_0 \circ \tilde{\mathscr{D}}\circ \pi_0&=& \lambda \mathscr{D} +\ad_{\varkappa};\\[2mm]
\text{\textup{(}Only when $p=2$\textup{)} } \mathscr{B}_\fg(x^*,x^*)&=&\lambda^{-2}( \mathscr{B}_\fa(\varkappa, \varkappa)+\mathscr{B}_{\tilde\fg}(\tilde x^*, \tilde x^*));
\end{eqnarray}
and 
\begin{equation}
\label{IsoEx}
\begin{array}{rcl}
\pi&=& \pi_0+ \mathscr{B}_\fa(\varkappa,\cdot) \tilde x \quad \text{on $\fa$;}\\[2mm]
\pi(x)&=&\lambda \tilde x;\\[2mm]
\pi(x^*)&=&\lambda^{-1}(\tilde x^*-\pi_0(\varkappa)-\rho \tilde x),\; \text{where}\\[3mm]
\rho&=&\left\{
\begin{array}{lcl}
\text{arbitrary} & \text{ if } & p=2,\\[2mm]
\frac{1}{2}\mathscr{B}_\fa(\varkappa,\varkappa) & \text{ if } & p\not =2.\\
\end{array}
\right.
\end{array}
\end{equation}
\end{Theorem}

\parbegin{Remark} {\rm 
If $\pi_0=\Id_\fa$, then condition (\ref{Cb}) means that $\mathscr{D}\simeq \tilde {\mathscr{D}}$ in $\mathrm{H}^1(\fa;\fa)$. Moreover, if the  derivations are restricted, then condition (\ref{Cb}) means $ {\mathscr{D}} \simeq \tilde {\mathscr{D}}$ in $\mathrm{H}^1_{\mathrm{res}}(\fa;\fa)$.}
\end{Remark}
Suppose now that $\fa$ is restricted with a~$p$-structure $[p]_\fa$. In Theorem \ref{MainTh} we proved that it is possible to extend the $p$-structure to any double extension. Let us denote by $[p]_\fg$ (resp. $[p]_{\tilde \fg}$) the $p$-structure on $\fg$ (resp. $\tilde \fg$) written in terms of $m,l,a_0, b_0, \gamma$ and $\mathscr{P}$ (resp. $\tilde m, \tilde l, \tilde a_0, \tilde b_0, \tilde \gamma$ and $\tilde{ \mathscr{P}}$). The following theorem characterizes the equivalence  class of $p$-structures on double extensions; we prove it only for $p=2,3$. 

An analog of Theorem \ref{pIsom1} for $p>3$ is still out of reach.

\parbegin{Theorem} [Isomorphism of DEs and $p$-homomorphism]\label{pIsom1}
For $p=2,3$, the  adapted isomorphism $\pi:  \fg\rightarrow  \tilde \fg$ given in Theorem \textup{\ref{Isom1}} defines a~$p$-homomorphism if and only if 
\[
\pi_0(a^{[p]_\fa})=(\pi_0(a))^{[p]_{\fa}}+\mathscr{B}_\fa(\varkappa,a)^p \tilde b_0,
\]
and
\[
\begin{array}{rcl}
\tilde m&=&\lambda^{-p} (\lambda m+\mathscr{B}_\fa(\varkappa,b_0)),\\[2mm]
\tilde b_0&=&\lambda^{-p} \pi_0(b_0),\\[2mm]
\tilde{ \mathscr{P}} \circ \pi_0&=&\lambda\mathscr{P}+\mathscr{B}_\fa(\varkappa, (\cdot)^{[p]_\fa})-\mathscr{B}_\fa(\varkappa,\cdot)^p\tilde m \\[2mm]
\tilde \gamma&=&\lambda^{p-1} \gamma,\\[2mm]
\tilde l&=&\lambda^{p}(\mathscr{B}_\fa(\varkappa, a_0)+\lambda l - \lambda^{-1} \gamma \rho) +\rho^p  \tilde m+\mathscr{P}(\pi_0(\varkappa))-\mathscr{B}_\fa(\tilde {\mathscr{D}}^{p-1}(\pi_0(\varkappa)),\pi_0(\varkappa)),\\[2mm]
\tilde a_0&=&\lambda^{p} (\pi_0(a_0)-\lambda^{-1}\gamma\pi_0(\varkappa))+(\pi_0(\varkappa))^{[p]_\fa}+\rho^p \tilde b_0+\tilde {\mathscr{D}}^{p-1}(\pi_0(\varkappa))\\[2mm]
&&-\left \{
\begin{array}{lll}
\lambda \pi_0([\mathscr{D}(\varkappa),\varkappa]_\fa),& \text{if} & $p=$3$$,\\[2mm]
0,&\text{if} & $p=$2$$.
\end{array}
\right.
\end{array}
\]
Moreover, if $\fz(\fa)=0$, then the automorphism $\pi_0$ of $\fa$ is also a~$p$-homomorphism.
\end{Theorem}
\begin{proof}
Let us study the conditions for which $\pi$ is a~$p$-homomorphism. We have
\[
\pi(x^{[p]_\fg})-(\pi(x))^{[p]_{\tilde \fg}}=\pi(m x+b_0)-(\lambda \tilde x)^{[p]_{\tilde \fg}}=m \lambda \tilde x+\pi_{0}(b_0)+\mathscr{B}_\fa(\varkappa,b_0)\tilde x-\lambda^p (\tilde m \tilde x+\tilde b_0).
\]
Therefore, $\tilde m=\lambda^{-p} (\lambda m+\mathscr{B}_\fa(\varkappa,b_0))$, and $\tilde b_0=\lambda^{-p} \pi_0(b_0)$.

Let $a\in \fa$. We have 
\[
\pi(a^{[p]_\fg})=\pi(a^{[p]_\fa}+\mathscr{P}(a)x)=\pi_0(a^{[p]_\fa})+\mathscr{B}_\fa(\varkappa, a^{[p]_\fa})\tilde x+\mathscr{P}(a)\lambda \tilde x.
\]
On the other hand, 
\[
(\pi(a))^{[p]_{\tilde \fg}}=(\pi_0(a)+\mathscr{B}_\fa(\varkappa,a)\tilde x)^{[p]_{\tilde \fg}}=(\pi_0(a))^{[p]_{\fa}}+\tilde{\mathscr{P}}(\pi_0(a))\tilde x+\mathscr{B}_\fa(\varkappa,a)^p (\tilde m \tilde x+\tilde b_0).
\]
It follows that (for every $a\in \fa$)
\begin{equation*}\label{pi01}
\begin{array}{l}
 \mathscr{B}_\fa(\varkappa,a)^p\tilde m +\tilde{\mathscr{P}}(\pi_0(a))-\mathscr{B}_\fa(\varkappa, a^{[p]_\fa})-\mathscr{P}(a)\lambda=0, \\[2mm] 
  \pi_0(a^{[p]_\fa})-(\pi_0(a))^{[p]_{\fa}}-\mathscr{B}_\fa(\varkappa,a)^p \tilde b_0=0.
\end{array}
\end{equation*}
Moreover,
\[
\begin{array}{lcl}
\pi((x^*)^{[p]_\fg})&=&\pi(a_0+l x+\gamma x^*)\\[2mm]
&=&\pi_0(a_0)+\mathscr{B}_\fa(\varkappa, a_0)\tilde x+l \lambda \tilde x+\gamma \lambda^{-1}(\tilde x^*-\pi_0(\varkappa)-\rho \tilde x)\\[2mm]
&=&(\mathscr{B}_\fa(\varkappa, a_0)+\lambda l - \lambda^{-1} \gamma \rho) \tilde x-\lambda^{-1} \gamma\pi_0(\varkappa)+\pi_0(a_0)+\lambda^{-1} \gamma \tilde x^*.
\end{array}
\]
On the other hand,
\[
\begin{array}{lcl}
(\pi(x^*))^{[p]_{\tilde \fg}}&=&\lambda^{-p}(\tilde x^*-\pi_0(\varkappa)- \rho \tilde x)^{[p]_{\tilde g}}\\[2mm]
&=&
\lambda^{-p} (\tilde a_0+\tilde \gamma \tilde x^*-\pi_0(\varkappa)^{[p]_\fa}-\rho^p \tilde b_0- \tilde {\mathscr{D}}^{p-1}(\pi_0(\varkappa))
\\[2mm]
&&+(\tilde l-\rho^p \tilde m-\tilde{\mathscr{P}}(\pi_0(\varkappa)))\tilde x) +\lambda^{-p}\mathscr{B}_\fa(\tilde {\mathscr{D}}^{p-1}(\pi_0(\varkappa)),\pi_0(\varkappa))\tilde x\\[2mm]
&&+ \left \{
\begin{array}{ll}
  \lambda^{-p} [\tilde {\mathscr{D}}(\pi_0(\varkappa)),\pi_0(\varkappa)]_\fa & \text{if $p=3$},\\[2mm]
0& \text{if $p=2$}.
\end{array}
\right .
\end{array}
\]
Therefore, 

\[
\begin{array}{rcl}
\lambda^{-1} \gamma&=&\lambda^{-p}\tilde \gamma,\\[2mm]
\mathscr{B}_\fa(\varkappa, a_0)+\lambda l - \lambda^{-1} \gamma \rho&=&\lambda^{-p}(\tilde l-\rho^p \tilde m-\mathscr{P}(\pi_0(\varkappa))+\mathscr{B}_\fa(\tilde {\mathscr{D}}^{p-1}(\pi_0(\varkappa)),\pi_0(\varkappa))),\\[2mm]
\pi_0(a_0)-\lambda^{-1}\gamma\pi_0(\varkappa)&=&\lambda^{-p} (\tilde a_0-\pi_0(\varkappa)^{[p]_\fa}-\rho^p \tilde b_0-\tilde {\mathscr{D}}^{p-1}(\pi_0(\varkappa)))\\[2mm]
&&+\left \{
\begin{array}{ll}
 \lambda^{-p} [\tilde {\mathscr{D}}(\pi_0(\varkappa)),\pi_0(\varkappa)]_\fa& \text{if $p=3$},\\[2mm]
0&\text{if $p=2$}.
\end{array}
\right.
\end{array}
\]
If $\fz(\fa)=0$, we have $b_0=\tilde b_0=0$ and the automorphism $\pi_0$ is a~$p$-homomorphism. 
\end{proof}

\subsubsection{The case of Lie superalgebras} \label{subsubsec:SIsom}
The main goal of this subsection is to superize Theorem \ref{Isom1}. The case $p=2$ has already been studied in \cite{BeBou}, so we assume that $p \not =2$. For $p=3$, the factor 3 in all expressions should be understood as zero.

Let $\pr: \tilde{\mathscr{K}} \oplus \fa\rightarrow \fa$ be the projection, and  ${\pi_0:=\pr \circ \pi}$. The map $\pi_0$ is obviously linear. Let $a\in \fa$. Since $\pi(a)-\pi_0(a)\in\mathrm{Ker}(\pr)$, it follows that $\pi(a)-\pi_0(a)\in  \tilde{\mathscr{K}} $. Since $\mathscr{B}_\fa$ is non-degenerate, there exists a~unique $t_\pi\in\fa_{p(\tilde{\mathscr{D}})}$ (depending only in $\pi$ and $p(\tilde{\mathscr{D}})$) such that 
\[
\pi(a)-\pi_0(a)=\mathscr{B}_\fa(t_\pi,a) \tilde x \quad \text{for any $a\in\fa$.}
\]

Besides, $\pi(x)=\lambda \tilde x$ for some $\lambda$ in $\mathbb{K}$. Indeed, let us write $\pi(x)=\lambda \tilde x+a$, where $a\in \fa$. We have 
\[
0=\mathscr{B}_\fg(x,b)=\mathscr{B}_{\tilde \fg}(\pi(x), \pi(b))=\mathscr{B}_{\tilde \fg}(\lambda \tilde x+a, \pi_0(b)+\mathscr{B}_\fa(t_\pi, b) \tilde x)=\mathscr{B}_\fa(a,\pi_0(b)). 
\]
Since $\pi_0$ is surjective and $\mathscr{B}_\fa$ is nondegenerate, it follows that $a=0$.
Let us show that $\pi_0$ preserves $\mathscr{B}_\fa$. Indeed, for any $a, b \in \fa$, we have
\[
\begin{array}{lcl}
\mathscr{B}_\fa(a,b)&=&\mathscr{B}_{\tilde \fg}(\pi(a), \pi(b))\\[2mm]
&=& \mathscr{B}_{\tilde \fg}(\pi_0(a)+\mathscr{B}_\fa(t_\pi, a)\tilde x, \pi_0(b)+\mathscr{B}_\fa(t_\pi, b)\tilde x)\\[2mm]
&=&\mathscr{B}_{\tilde \fg}(\pi_0(a), \pi_0(b)) 
=\mathscr{B}_{\fa}(\pi_0(a), \pi_0(b)).
\end{array}
\] 
Let us show that $\pi_0$ is an automorphism of $(\fa, \mathscr{B}_\fa)$.  Let $a, b \in \fa$. We get 
\[
\begin{array}{lcl}
\pi([a,b]_\fg)&=&\pi([a,b]_\fa)+(-1)^{p(\mathscr{D})(p(\mathscr{B}_\fa)+1)}\lambda \mathscr{B}_\fa(\mathscr{D}(a),b)\tilde x\\[2mm]
&=&\pi_0([a,b]_\fa)+\mathscr{B}_\fa(t_\pi,[a,b]_\fa)\tilde x+(-1)^{p(\mathscr{D})(p(\mathscr{B}_\fa)+1)}\lambda \mathscr{B}_\fa(\mathscr{D}(a),b)\tilde x
\end{array}
\]
and 
\[
\begin{array}{lcl}
[\pi(a), \pi(b)]_{\tilde \fg}&=&[\pi_0(a)+\mathscr{B}_\fa(t_\pi,a)\tilde x, \pi_0(b)+\mathscr{B}_\fa(t_\pi,b) \tilde x]_{\tilde \fg}\\[2mm]
&=&[\pi_0(a), \pi_0(b)]_{\tilde \fg} \\[2mm]
&=& [\pi_0(a), \pi_0(b)]_{\fa}+(-1)^{p(\tilde{\mathscr{D}})(p(\mathscr{B}_\fa)+1)}\mathscr{B}_\fa(\tilde{\mathscr{D}}(\pi_0(a)),\pi_0(b))\tilde x.
\end{array}
\]
It follows that 
\[
\mathscr{B}_\fa(\pi_0^{-1}\tilde{\mathscr{D}}(\pi_0(a)),b)=(-1)^{p(\mathscr{D})(p(\mathscr{B}_\fa)+1)}\mathscr{B}_\fa([t_\pi,a]_\fa,b)+\lambda \mathscr{B}_\fa(\mathscr{D}(a),b).
\] 
Therefore, 
\begin{equation}
\label{DDt}
\pi_0^{-1}\circ \tilde{\mathscr{D}} \circ \pi_0 =(-1)^{p(\mathscr{D})(p(\mathscr{B}_\fa)+1)} \ad_{t_\pi} +\lambda \mathscr{D}.
\end{equation}

Let us write 
\[
\pi(x^*)= \mu \tilde x^*+ a+ 
 \nu \tilde x, \quad \text{where $a \in \fa$ and $\nu=0$ if $\mathscr{B}_\fa$ is odd.
}
\]
Because $\pi$ is an isomorphism, we have
\[
1=\mathscr{B}_{\tilde \fg}(\pi(x), \pi(x^*))=\mathscr{B}_{\tilde \fg}( \lambda \tilde x, \mu \tilde x^*+a+\nu \tilde x)=\lambda \mu.
\]
Therefore, $\mu=\lambda^{-1}$. Besides,  
\[
\begin{array}{lcl}
0=\mathscr{B}_{\tilde \fg}(\pi(x^*), \pi(b))&=&\mathscr{B}_{\tilde \fg}(\tilde x^*+a+\nu \tilde x,  \pi_0(b)+\mathscr{B}_\fa(t_\pi,b)\tilde x)\\[2mm]
&=&\mu (-1)^{p(\tilde x)p(\tilde x^*)}\mathscr{B}_\fa(t_\pi,b)+\mathscr{B}_\fa(a, \pi_0(b)).
\end{array}
\]
Hence, $ \mathscr{B}_\fa(\mu  (-1)^{p(\tilde x)p(\tilde x^*)} t_\pi+\pi_0^{-1}a, b)=0$, and $\mu  (-1)^{p(\tilde x)p(\tilde x^*)} t_\pi+\pi_0^{-1}(a)=0$ implying 
\begin{equation}
a= -(-1)^{p(\tilde x)p(\tilde x^*)}\mu \pi_0(t_\pi).
\end{equation} 
Finally, we get
\[
\pi(x^*)=\lambda^{-1}(\tilde x^*- (-1)^{p(\tilde x)p(\tilde x^*)}\pi_0(t_\pi))+\nu \tilde x.
\]
Besides,
\[
0=\mathscr{B}_{\tilde \fg}(\pi(x^*), \pi(x^*))=\mathscr{B}_{\tilde \fg}( \mu \tilde x^*+a+\nu \tilde x, \mu \tilde x^*+a+\nu \tilde x)=\mu \nu(1+(-1)^{p(\tilde x) p(\tilde x^*)} )+\mathscr{B}_\fa(a,a).
\]
We distinguish two cases: 

(i) If $\mathscr{B}_\fa$ is odd, then $\nu=0$, and therefore $\mathscr{B}_\fa(a,a)=0$. 

(ii) If $\mathscr{B}_\fa$ is even, then $\mathscr{B}_\fa(a,a)=0$ and $\nu$ is arbitrary, except when $p(x)=p(x^*)=\ev$. In the latter case (recall that $p\not =2$)
\[
\nu=-\frac{1}{2}\lambda^{-1} \mathscr{B}_\fa(t_\pi, t_\pi).
\] 

Let $p(x^*)=\od$ (and therefore $p(\mathscr{D})=p(\tilde{\mathscr{D}})=\od$). 
We have
\[
\begin{array}{ll}
[\pi(x^*), \pi(x^*)]_{\tilde \fg}&=[\mu \tilde x^*+a+\nu x, \ \mu \tilde x^*+a+\nu \tilde x]_{\tilde \fg}\\[2mm]
&=2 \mu ( \tilde{\mathscr{D}}(a)-(-1)^{p(a)p(\mathscr{B}_\fa)}2 \mathscr{B}_\fa(a, \tilde b_0)\tilde x)+[a,a]_\fa\\
&-(-1)^{p(\mathscr{B}_\fa)} \mathscr{B}_\fa(\tilde{\mathscr{D}}(a),a) \tilde x+\mu^2 [\tilde x^*, \tilde x^*]_{\tilde \fg}.
\end{array}
\]

On the other hand,
\[
\pi([x^*, x^*]_\fg)=\left \{
\begin{array}{ll}
2\pi_0(b_0)+2\mathscr{B}_\fa(t_\pi, b_0)\tilde x+\lambda \lambda_0 \tilde x & \text{ if  $p(\mathscr{B}_\fa)=\od$,}\\[2mm]
2\pi_0(b_0) & \text{ if $p(\mathscr{B}_\fa)=\ev$.}
\end{array}
\right. 
\]
It follows that 
\begin{equation} 
2\tilde b_0= \mu^{-2}(2 \pi_0(b_0)-[a,a]_\fa-2\mu \tilde{\mathscr{D}}(a)).
\end{equation} 
Therefore, 
\begin{equation} 
\tilde b_0= \lambda^2 \pi_0(b_0)+ (-1)^{p(x)} \lambda \pi_0(\mathscr{D}(t_\pi))+ \frac{1}{2}\pi_0([t_\pi,t_\pi]_\fa).
\end{equation} 
In addition, in the case where $p(\mathscr{B}_\fa)=\od$ we have (see, Eq. (\ref{lambda}))
\[
\lambda_0 \lambda +2 \mathscr{B}_\fa(t_\pi,b_0)=\mu^2 \tilde \lambda_0+4 \mu \mathscr{B}_\fa(a,\tilde b_0)+\mathscr{B}_\fa(\tilde{\mathscr{D}}(a),a). 
\]
After computation, we get
\[
\tilde \lambda_0=\lambda_0 \lambda^3 +\mathscr{B}_\fa(t_\pi, [t_\pi,t_\pi]_\fa)+3 \lambda \mathscr{B}_\fa (t_\pi, \mathscr{D}(t_\pi))+ 6 \lambda^2 \mathscr{B}_\fa(t_\pi, b_0).
\]
We arrive at the following two theorems.

\parbegin{Theorem}[Isomorphism of DEs for even derivations] \label{SIsom2}
Let $\mathscr{B}_\fa$ be $\mathscr{D}$- and $\tilde {\mathscr{D}}$-invariant, where $\mathscr{D}, \tilde {\mathscr{D}} \in \fder_\ev(\fa)$ (and satisfy condition \textup{(\ref{p=3Con2})} if $p=3$ and $p(\mathscr{B}_\fa)=\od$).

Then, there exists an adapted isomorphism $\pi:  \fg\rightarrow  \tilde \fg$ if and only if there exists an automorphism $\pi_0:\fa\rightarrow \fa$, a~$\lambda \in \mathbb{K} ^{\times}$, and $\varkappa \in \fa$ \textup{(}all unique\textup{)} such that
\[
\begin{array}{rcl}
\pi_0^{-1}\circ \tilde{\mathscr{D}} \circ \pi_0&=&\lambda \mathscr{D} + \ad_{\varkappa};\\[2mm]
\mathscr{B}_\fa(\varkappa, \varkappa)&=&0\quad \text{if \; $p(\mathscr{B}_\fa)=\od$}.\\[2mm]
\end{array}
\]
and  
\[
\begin{array}{rcl}
\pi&=& \pi_0+ \mathscr{B}_\fa(\varkappa,\cdot) \tilde x \quad \text{on $\fa$}; \\[2mm]
\pi(x)&=&\lambda \tilde x;\\[2mm]
\pi(x^*)&=&\lambda^{-1}(\tilde x^*- \pi_0(\varkappa)-\frac{1}{4}(1+(-1)^{p(\mathscr{B}_\fa)} )  \mathscr{B}_\fa(\varkappa, \varkappa)x ).
\end{array}
\]
\end{Theorem}

\parbegin{Theorem}[Isomorphism of DEs for odd derivations] \label{SIsom3} 
Let $\mathscr{B}_\fa$ be $\mathscr{D}$- and $\tilde {\mathscr{D}}$-invariant, where $\mathscr{D}, \tilde {\mathscr{D}} \in \fder_\od(\fa)$ (and satisfy condition \textup{(\ref{p=3Con1})} if $p=3$ and $p(\mathscr{B}_\fa)=\od$). Then, there exists an adapted isomorphism $\pi:  \fg\rightarrow  \tilde \fg$ if and only if there exists an automorphism $\pi_0:\fa\rightarrow \fa$, a~ $\lambda \in \mathbb{K} ^{\times}$, and $\varkappa \in \fa$ \textup{(}all unique\textup{)} such that
\[
\begin{array}{rcl}
\pi_0^{-1} \circ \tilde{\mathscr{D}} \circ \pi_0&=&\lambda \mathscr{D}-(-1)^{p(\mathscr{B}_\fa)}\ad_{\varkappa};\\[2mm]
\mathscr{B}_\fa(\varkappa,\varkappa)&=&0;\\[2mm]
\end{array}
\]
and 
\[
\begin{array}{rcl}
\pi&=& \pi_0+ \mathscr{B}_\fa(\varkappa,\cdot) \tilde x \quad \text{on $\fa$}; \\[2mm]
\pi(x)&=&\lambda \tilde x;\\[2mm]
\pi(x^*)&=&\lambda^{-1}(\tilde x^*- (-1)^{p(x) }\pi_0(\varkappa)) + 
(1+(-1)^{p(\mathscr{B}_\fa)} )  \nu x, 
\text{where $\nu$ is arbitrary}.
\end{array}
\]
Moreover, 
\[
\begin{array}{lcl}
\tilde b_0&=& \lambda^2 \pi_0(b_0) +(-1)^{p(x)}  \lambda (\pi_0(\mathscr{D}(\varkappa)))+ \frac{1}{2} \pi_0([\varkappa,\varkappa]_\fa).
\end{array}
\]

If $p(\mathscr{B}_\fa)=\od$, we additionally have 
 \[
\begin{array}{lcl}
\tilde \lambda_0=\lambda_0 \lambda^3 +\mathscr{B}_\fa(\varkappa, [\varkappa,\varkappa])+3 \lambda \mathscr{B}_\fa (\varkappa, \mathscr{D}(\varkappa))+ 6 \lambda^2 \mathscr{B}_\fa(\varkappa, b_0).
\end{array}
\]
\end{Theorem}

\begin{proof}[Proof of Theorem $\ref{SIsom2}$ and Theorem $\ref{SIsom3}$] To check that $\pi$ preserves the Lie bracket, it remains to check that $\pi([x^*,a]_\fg)=[\pi(x^*),\pi(a)]_{\tilde \fg}$. We distinguish two cases:

\underline{(i) The case where $p(x^*)=\od$.}  We have 
\[
\begin{array}{lcl}
\pi([x^*,a]_\fg)&=&\pi (\mathscr{D}(a) -(-1)^{p(\mathscr{B}_\fa) p(a)} 2 \mathscr{B}_\fa(a,b_0) x)\\
&=&\pi_0(\mathscr{D}(a))+ \mathscr{B}_\fa(\varkappa, \mathscr{D}(a))\tilde x 
-(-1)^{p(\mathscr{B}_\fa) p(a)} 2 \lambda \mathscr{B}_\fa(a,b_0) \tilde x.
\end{array}
\] 
On the other hand, 


\[
\begin{array}{lcl}
[\pi(x^*),\pi(a)]_{\tilde \fg}&=&
\pi_0 (\mathscr{D}(a))-(-1)^{p(a)p(\mathscr{B}_\fa)} 2 \lambda \mathscr{B}_\fa(a, b_0) \tilde x \\[2mm]
&& +(-2 (-1)^{p(x)+p(a)p(\mathscr{B}_\fa)}-1)\mathscr{B}_\fa( a, \mathscr{D}(\varkappa)) \tilde x\\[2mm]
&&+\lambda^{-1}((-1)^{p(\mathscr{B}_\fa)}-(-1)^{p(a)p(\mathscr{B}_\fa)}) \mathscr{B}_\fa(a,[\varkappa,\varkappa]_\fa)\tilde x.

\end{array}
\]
The result follows since $p(\mathscr{B}_\fa)=p(x)+p(x^*)=p(x)+1$, and (for all $a\in \fa$)
\[
((-1)^{p(\mathscr{B}_\fa)}-(-1)^{p(a)p(\mathscr{B}_\fa)}) \mathscr{B}_\fa(a,[\varkappa,\varkappa]_\fa) =0.
\] 
\underline{(ii) The case where $p(x^*)=\ev$.} Similar computation.

There is no need to check the remaining brackets because they are certainly satisfied as shown by the previous computations prior formulations of Theorems \ref{SIsom2} and  \ref{SIsom3}.

Let us show that $\pi$ preserves $\mathscr{B}_\fg$. For every element $a\in \fa$, we have $\mathscr{B}_\fg(x^*,a)=0$. On the other hand,
\[
\begin{array}{lcl}
\mathscr{B}_{\tilde \fg}(\pi(a), \pi(x^*))&=&\mathscr{B}_{\tilde \fg}(\pi_0(a)+\mathscr{B}_\fa(\varkappa,a )\tilde x, \pi(x^*))\\[2mm]
&=&\lambda^{-1}\mathscr{B}_\fa(\varkappa, a)-(-1)^{p(x)p(x^*)}\mathscr{B}_{\fa} (\pi_0(a), \lambda^{-1}\pi_0(\varkappa))\\[2mm]
&=&\lambda^{-1} ((-1)^{p(x)p(a)}-(-1)^{p(x)p(x^*)})\mathscr{B}_\fa(\varkappa, a) =0.
\end{array}
\]
Let $\mathscr{D}(b_0)=0$. Let us show that $\tilde{\mathscr{D}}(\tilde b_0)=0$. Indeed, (recall that $p(\mathscr{D})=p(\tilde{\mathscr{D}})=\od$): 
\[
\begin{array}{lcl}
\pi_0^{-1}\tilde{\mathscr{D}}(\tilde b_0)&=&\pi_0^{-1}\left ( \tilde{\mathscr{D}} ( \lambda^2 \pi_0(b_0) +(-1)^{p(x)p(x^*)} \lambda (\pi_0(\mathscr{D}(\varkappa)))+ \frac{1}{2}\pi_0([\varkappa,\varkappa]_\fa)\right )\\[2mm]
&=&(\lambda \mathscr{D}-(-1)^{p(\mathscr{B}_\fa)}\ad_{\varkappa}) \left ( \lambda^2 b_0 +(-1)^{p(x)p(x^*)} \lambda \mathscr{D}(\varkappa)+ \frac{1}{2} [\varkappa,\varkappa]_\fa\right ) =0.
\end{array}
\]
Let us show that $\tilde{\mathscr{D}}^2=\ad_{\tilde b_0}$. Indeed,
\[
\begin{array}{lcl}
\tilde{\mathscr{D}}^2&=& \pi_0 (\lambda \mathscr{D} - (-1)^{p(\mathscr{B}_\fa)}\ad_\varkappa)^2\pi_0^{-1}\\[2mm]
&=& \pi_0 (\lambda^2 \mathscr{D}^2 - (-1)^{p(\mathscr{B}_\fa)}( \lambda \mathscr{D} \ad_\varkappa+ \lambda \ad_\varkappa \mathscr{D})+\ad_\varkappa^2)\pi_0^{-1}\\[2mm]
&=& \pi_0 (\lambda^2 \ad_{b_0}- \lambda(-1)^{p(\mathscr{B}_\fa)}  \ad_{\mathscr{D}(\varkappa)}+ \frac{1}{2} \ad_{[\varkappa, \varkappa]_\fa})\pi_0^{-1}= \ad_{\tilde b_0}.
\end{array}
\]

Besides, (since $\tilde b_0=\lambda^2 \pi_0(b_0) +(-1)^{p(x)}  \lambda (\pi_0(\mathscr{D}(\varkappa)))+  \frac{1}{2}\pi_0([\varkappa,\varkappa]_\fa)$ and $\pi_0$ is an isomorphism):
\[
\begin{array}{lcl}
\mathscr{B}_{\tilde \fg}(\tilde b_0, \tilde b_0)&=&
 \lambda^4 \mathscr{B}_\fa(b_0,b_0)+  \lambda^2 \mathscr{B}_\fa(b_0, [\varkappa,\varkappa]_\fa)+  \lambda^2 \mathscr{B}_\fa(\mathscr{D}(\varkappa),\mathscr{D}(\varkappa))\\[2mm]
 &&+(-1)^{p(x)} \lambda \mathscr{B}_\fa (\mathscr{D}(\varkappa), [\varkappa,\varkappa]_\fa)=0.
 \qed
\end{array}
\]
\phantom\qedhere 
\end{proof}

Suppose now that $\fa$ is restricted. In Theorems \ref{DoBe} and \ref{DoBo}, we proved that it is possible to extend its $p|2p$-structure to any of its double extension. Let us denote by $[p|2p]_\fg$ (resp. $[p|2p]_{\tilde \fg}$) the $p|2p$-structure on $\fg$ (resp. $\tilde \fg$) written in terms of $m,l,a_0, b_0, \gamma$ and $\mathscr{P}$ (resp. $\tilde m, \tilde l, \tilde a_0, \tilde b_0, \tilde{ \gamma}$ and $\tilde{\mathscr{P}}$). The following theorems characterize the equivalence  class of $p|2p$-structures on double extensions. Theorem \ref{pIsom2} is proved only for $p=3$. 

\parbegin{Theorem}[Isomorphism of DEs and $p$-homomorphism for even derivations] \label{pIsom2} 
The adapted isomorphism $\pi:  \fg\rightarrow  \tilde \fg$ given in Theorem \textup{\ref{SIsom2}} defines a~$p|2p$-homomorphism, if and only if the equations below are satisfied: 

\textup{(}i\textup{)} If $p(\mathscr{B}_\fa)=\ev$ (here $p=3$), we have 
\[
\pi_0(a^{[p]_\fa})=(\pi_0(a))^{[p]_{\fa}}+  \mathscr{B}_\fa(\varkappa,a)^p \tilde c_0 \quad \text{for all $a\in \fa_\ev$,}
\]
and
\[
\begin{array}{rcl}
\tilde m&=&\lambda^{-p} (\lambda m+ \mathscr{B}_\fa(\varkappa,c_0)),\\[2mm]
\tilde c_0&=&\lambda^{-p} \pi_0(c_0),\\[2mm]
\tilde{\mathscr{P}} \circ \pi_0&=&\lambda\mathscr{P}+\mathscr{B}_\fa(\varkappa, (\cdot)^{[p]_\fa})-\mathscr{B}_\fa(\varkappa,\cdot)^p\tilde m \text{ on $\fa_\ev$}\\[2mm]
\tilde \gamma&=&\lambda^{p-1} \gamma,\\[2mm]
\tilde l&=&\lambda^{p}(\mathscr{B}_\fa(\varkappa, a_0)+\lambda l - \lambda^{-1} \gamma \rho) +\rho^p  \tilde m+\mathscr{P}(\pi_0(\varkappa))-\mathscr{B}_\fa(\tilde{\mathscr{D}}^{p-1}(\pi_0(\varkappa)),\pi_0(\varkappa))),\\[2mm]
\tilde a_0&=&\lambda^{p} (\pi_0(a_0)-\lambda^{-1}\gamma\pi_0(\varkappa))+(\pi_0(\varkappa))^{[p]_\fa}+  \rho^p \tilde c_0+\tilde{\mathscr{D}}^{p-1}(\pi_0(\varkappa))- \lambda \pi_0([\mathscr{D}(\varkappa),\varkappa]_\fa).
\end{array}
\]
Moreover, if $\fz(\fa)=0$, then the automorphism $\pi_0$ of $\fa$ is a~also a~$p|2p$-homomorphism.

\textup{(}ii\textup{)} If $p(\mathscr{B}_\fa)=\od$ (here $p=3$), we have
\[
\pi_0(a^{[p]_\fa})=(\pi_0(a))^{[p]_{\fa}} \quad \text{for all $a\in \fa_\ev$,}
\]
and
\[
\begin{array}{rcl}
\tilde \gamma&=&\lambda^{p-1} \gamma ,\\[3mm]
\tilde a_0&=& \pi_0(\varkappa^{[p]_\fa})+\lambda^2 \pi_0 (\mathscr{D}^2(\varkappa))-\gamma \lambda^{p-1} \pi_0(\varkappa) + \lambda^p\pi_0(a_0)+ 2 \lambda \pi_0([\varkappa,\mathscr{D}(\varkappa)]_\fa).
\end{array}
\]

\end{Theorem}

\parbegin{Theorem}[Isomorphism of DEs and $p$-homomorphism for odd derivations] \label{pIsom3} 
The  adapted isomorphism $\pi:  \fg\rightarrow  \tilde \fg$ given in Theorem \textup{\ref{SIsom3}} defines a~$p|2p$-homomorphism, if and only if the following conditions are satisfied:

\textup{(}i\textup{)} If $p(\mathscr{B}_\fa)=\od$, then 

\[
\pi_0(a^{[p]_\fa})=(\pi_0(a))^{[p]_{\fa}} +  \mathscr{B}_\fa(\varkappa,a)^p \tilde c_0\quad \text{for all $a\in \fa_\ev$},
\]
and
\[
\begin{array}{rcl}
\tilde m&=&\lambda^{-p} (\lambda m+ \mathscr{B}_\fa(\varkappa,c_0)),\\[2mm]
\tilde c_0&=&\lambda^{-p} \pi_0(c_0),\\[2mm]
\tilde{\mathscr{P}} \circ \pi_0&=&\lambda\mathscr{P}+\mathscr{B}_\fa(\varkappa, (\cdot)^{[p]_\fa})-\mathscr{B}_\fa(\varkappa,\cdot)^p\tilde m\quad  \text{ on $\fa_\ev$}.\\ 
\end{array}
\]
Moreover, if $\fz(\fa)=0$, then the automorphism $\pi_0$ of $\fa$ is also a~$p|2p$-homomorphism.

\textup{(}ii\textup{)} If $p(\mathscr{B}_\fa)=\ev$, then $\pi_0$ defines a~$p|2p$-homomorphism on $\fa$.
\end{Theorem}
\begin{proof} The proof of Theorems  \ref{pIsom2} and \ref{pIsom3} is similar to that of Theorem \ref{pIsom1}, minding the Sign Rule.
\end{proof}
\section{Examples}\label{sec: exa}

{\rm The calculations performed with the aid of \textit{SuperLie} code \cite{Gr} are called \textbf{Claims}. The 1-cochain $\hat x\in C^1(\fg)$ denotes the dual of $x\in\fg$.}  

\subsection{$\fpsl(3)$ for $p=3$}
Let us fix a~basis of $\mathfrak{psl}(3)$ generated by the root vectors $x_1$, $x_2$, $x_3=[x_1,x_2]$ (positive) and $y_1$, $y_2$, $y_3=[y_1,y_2]$ (negative). In the ordered basis $e_1=[x_1,y_1]$, $e_2= x_1$, $e_3=x_2$, $e_4=x_3$, $e_5=y_1$, $e_6=y_2$, $e_7=y_3 $ of $\mathfrak{psl}(3)$ NIS has the Gram matrix
\[
\mathscr{B}_{\mathfrak{psl}(3)}=\left (
\begin{array}{ccc}
-1& 0 &0\\
0& 0 & I_{2,1}\\
0 &I_{2,1}&0
\end{array}
\right ),
\]
where $I_{r,s}=\text{diag}(1,\ldots,1,-1,\ldots,-1)$ with $r$-many $1$s and $s$ many $(-1)$s, see \cite{BKLS}. The $3$-structure on $\mathfrak{psl}(3)$ is given by the formulas:
\[
e_1^{[3]}=e_1, \quad e_i^{[3]}=0 \quad \text{ for all $i>1$}.
\]
\sssbegin{Claim} [\cite{Ib2}] The  space $\mathrm{H}^1_{\mathrm{res}}(\fpsl(3); \fpsl(3))$ is spanned by the $7$ cocycles of which we only need the following ones, their degrees in the subscript. \textup{(In fact, we have $\mathfrak{out}(\mathfrak{psl}(3)) \simeq \mathfrak{psl}(3)$.)}
\begin{equation}\label{cocypsl3}\tiny
\begin{array}{l}
\mathscr{D}_{-3}^1=y_1 \otimes \widehat{x_3}+
   y_3 \otimes \widehat{x_1},\\[1mm]  
   
 
 \mathscr{D}_{0}^2=2\, x_1 \otimes \widehat{x_2}+ y_2 \otimes \widehat{y_1}, \\ [1mm]
 \mathscr{D}_{0}^3=x_1 \otimes \widehat{x_1}+
    x_3 \otimes \widehat{x_3}+
   2\, y_1 \otimes \widehat{y_1}+
   2\, y_3 \otimes \widehat{y_3},\\ 

\end{array}
\end{equation}
\end{Claim}

An easy computation shows that $\mathscr{B}_{\mathfrak{psl}(3)}$ is $\mathscr{D}$-invariant for any  $\mathscr{D}$  in \eqref{cocypsl3}.

The isomorphism  
\[
y_1 \leftrightarrow x_1, \; y_2 \leftrightarrow x_2, \;y_1 \leftrightarrow y_2, \; y_3 \leftrightarrow -y_3, \;
\]
sends $\mathscr{D}_{-3}^1$ to a new one $\mathscr{D}_{3}^1$. The root system of $\mathfrak{psl}(3)$
is symmetric, so the isomorphism
\[
\displaystyle x_1 \leftrightarrow x_2, \; x_3 \leftrightarrow -x_3, \;y_1 \leftrightarrow y_2, \; y_3 \leftrightarrow -y_3, \;
\]
sends the cocycles described in  \eqref{cocypsl3} and $\mathscr{D}_{-3}^1$ into the new ones:
\[
\mathscr{D}_{0}^2 \leftrightarrow  \mathscr{D}_{0}^1, \; \mathscr{D}_{-3}^1 \leftrightarrow \mathscr{D}_{-3}^2, \; \mathscr{D}_3^{1} \leftrightarrow \mathscr{D}_{3}^2.
\] 
Therefore, we can confine ourselves to $\mathscr{D}_{-3}^1$, $\mathscr{D}_{0}^2$ and $\mathscr{D}_{0}^3$ only. Let $\mathscr{P}_i$ be the corresponding cubic forms on $\fa$ defined as follows: 
\[
\begin{array}{cll}
\mathscr{P}_1  (\sum_{1\leq i \leq 7} \lambda_i e_i  )&:=&\lambda _6 \lambda _4^2+2 \lambda _1 \lambda _2 \lambda_4+2\lambda _3 \lambda _2^2,\\[2mm]
\mathscr{P}_2  (\sum_{1\leq i \leq 7} \lambda_i e_i  )&:=&\lambda _7 \lambda _3^2+2 \lambda _1 \lambda _5 \lambda
   _3+\lambda _4 \lambda _5^2,\\[2mm]
 \mathscr{P}_3  (\sum_{1\leq i \leq 7} \lambda_i e_i )&:=&\lambda _1 \lambda _2 \lambda _5+\lambda _4 \lambda _6
   \lambda _5+2 \lambda _2 \lambda _3 \lambda _7+\lambda
   _1 \lambda _4 \lambda _7.
\end{array}
\] 

Now, let $a=\sum_{1\leq i \leq 7} \lambda_i e_i,\; b=\sum_{1\leq i \leq 7} \mu_i e_i \in \mathfrak{psl}(3).$ We have 
\[
\begin{array}{lcl}
\mathscr{B}_{\mathfrak{psl}(3)}(\mathscr{D}_{3}^0(a-b), [a,b])&=&(\lambda_2-\mu_2)(\lambda_1\mu_5+\lambda_7\mu_3-\lambda_5\mu_1-\lambda_3\mu_7)\\[2mm]
&&-(\lambda_4-\mu_4)(\lambda_7\mu_1+\lambda_5\mu_6-\lambda_1\mu_7-\lambda_6\mu_5)\\[2mm]
&&-(\lambda_5-\mu_5)(\lambda_2\mu_1+\lambda_4\mu_6-\lambda_1\mu_2-\lambda_6\mu_4)\\[2mm]
&&+(\lambda_7-\mu_7)(\lambda_2\mu_3+\lambda_1\mu_4-\lambda_3\mu_2-\lambda_4\mu_1).
\end{array}
\]
Besides,
\[
\begin{array}{lcl}
\mathscr{B}_{\mathfrak{psl}(3)}(\mathscr{D}_0^2(a-b), [a,b])&=&-(\lambda_3-\mu_3)(\lambda_1\mu_5+\lambda_7\mu_3-\lambda_5\mu_1-\lambda_3\mu_7)\\[2mm]
&&+(\lambda_5-\mu_5)(\lambda_3\mu_1+\lambda_5\mu_4-\lambda_1\mu_3-\lambda_4\mu_5),\\[2mm]
\mathscr{B}_{\mathfrak{psl}(3)}(\mathscr{D}_{-3}^1(a-b), [a,b])&=&-(\lambda_2-\mu_2)(\lambda_1\mu_4+\lambda_2\mu_3-\lambda_3\mu_2-\lambda_4\mu_1)\\[2mm]
&&+(\lambda_4-\mu_4)(\lambda_2\mu_1+\lambda_4\mu_6-\lambda_1\mu_2-\lambda_6\mu_4).
\end{array}
\]
%
Let us summarize and introduce notation: 
\begin{equation}\label{psl3}
\footnotesize
\renewcommand{\arraystretch}{1.4}
\begin{tabular}{|c|c|c|c|c|} \hline
Derivation & $\mathscr{P}(a)$&$\gamma$ &$a_0$&
Double extension\\
\hline
$\mathscr{D}_0^2$ & $\lambda _3^2 \lambda _7 +2 \lambda _1 \lambda _5 \lambda
   _3+\lambda _4 \lambda _5^2$ &0&0& $\widetilde \fgl(3)$ \\ \hline
   $\mathscr{D}_{-3}^1$ & $\lambda_4^2\lambda_6+2 \lambda_1\lambda_2\lambda_4+2 \lambda_2^2\lambda_3$ &0&0& $ \widehat{\fgl}(3)$ \\ \hline
$\mathscr{D}_0^3$ & $\lambda _1 \lambda _2 \lambda _5+\lambda _4 \lambda _6
   \lambda _5+2 \lambda _2 \lambda _3 \lambda _7+\lambda
   _1 \lambda _4 \lambda _7$ &1& 0&$\fgl(3)$\\ \hline
\end{tabular}
\end{equation}

\sssbegin{Claim}  $\dim \mathrm{H}^2(\fgl(3))=0$, $\dim \mathrm{H}^2(\widetilde \fgl(3))=3$, and $\dim \mathrm{H}^2(\widehat{\fgl}(3))=4$, hence  $\fgl(3)$, $\widetilde{\fgl}(3)$, and $\widehat{\fgl}(3)$ are pairwise non-isomorphic. 
\end{Claim}
\subsection{$\fbrj(2;3)$ for $p=3$}\label{prebrj}
Let us realize $\mathfrak{brj}(2;3)$ (for more details, see \cite{BGL, BoGL}) by the Cartan matrix and the positive root vectors (odd ones are underlined)
\[
\left (
\begin{array}{rr}
 0 & -1\\
-2 & 1
\end{array}
\right ) \quad 
\begin{array}{l}
\underline{x_1}, \underline{x_2}, \  x_3=[x_1, x_2], \  x_4=[x_2, x_2] , \  \underline{x_5}  = [x_2,x_3],\\[2mm]
x_6= [x_3, x_4] , \ \ \underline{x_7}= [x_4, x_5], \ \ x_8= [x_5, x_5]. \\
\end{array}
\]
As shown in \cite{BKLS}, the Lie superalgebra $\mathfrak{brj}(2;3)$ admits a~NIS given in the ordered basis 
\begin{equation}
\label{basbrj23}
h_1:=[x_1,y_1], \quad h_2:=[x_2,y_2], \quad x_1,\ldots, x_8, \quad y_1,  \ldots, y_8. 
\end{equation}
by the Gram matrix
\[
\begin{array}{l}
\mathscr{B}_{\mathfrak{brj}(2;3)}=\left (
\begin{array}{ccc}
A& 0 &0   \\
0& 0  & B\\
0 & C &0 \\
\end{array}
\right ),\text{~~where $A=\left ( 
\begin{array}{cc} 
0 & 1\\
1 & 2
\end{array}\right )$,} 
\\
\text{$B=\mathrm{diag} \{1,2, 2, 2, 2, 2, 1,1\}$ and $C=\mathrm{diag} \{2, 1, 2, 2,1, 2,2,1\}$.}
\end{array}
\]

There exists a~$3|6$ structure on $\fbrj(2;3)$ that we express in the basis (\ref{basbrj23}) as follows: 
\[
h_1^{[3]}=h_1,\,h_2^{[3]}=h_2, \, v^{[3]}= w^{[6]}=0 \text{ for any root vectors $v$ (even) and $w$ (odd) in (\ref{basbrj23})}
\]
\sssbegin{Claim} The  space $\mathrm{H}^1 (\fbrj(2;3); \fbrj(2;3))$ is spanned by the odd cocycles\footnote{Recall that $c(a_k)=(-1)^{p(a_k)} \delta_j^k a_i$ for the 1-cochain $c=a_i\otimes \widehat{a_j}$, }: 
\begin{equation}\label{cocybrj3}\tiny
\begin{array}{ll}
\deg=-3:&\mathscr{D}_{-3}= x_1 \otimes  \widehat{x_6}+
   x_3 \otimes \widehat{x_7} +2 y_2 \otimes \widehat{x_4}+ y_4 \otimes \widehat{x_2}+ 2 y_6 \otimes \widehat{y_1}+ y_7 \otimes \widehat{y_3}, \\[2mm]

\deg=3: &\mathscr{D}_{3}=x_2  \otimes \widehat{y_4} + 
 x_4 \otimes  \widehat{y_2} + 
 x_6 \otimes \widehat{x_1} + 
 2 x_7 \otimes \widehat{x_3} + 
 y_1 \otimes \widehat{y_6} + 
y_3 \otimes \widehat{y_7}.
\end{array}
\end{equation}
\end{Claim}
A direct computation shows that $\mathscr{B}_{\fbrj (2;3)}$ is $\mathscr{D}$-invariant, where $\mathscr{D}$ is any of the derivations given in (\ref{cocybrj3}). However, the condition (\ref{p=3Con1}) is violated for these derivations. Indeed, 
\[
\begin{array}{lcl}
\mathscr{B}_\fa(\mathscr{D}_{-3}(x_2),[x_2,x_2])&=&\mathscr{B}_\fa(2y_4, x_4)=1\not=0, \\[2mm] \mathscr{B}_\fa(\mathscr{D}_{3}(y_2),[y_2,y_2])&=&\mathscr{B}_\fa(2x_4, y_4)=1\not=0.
\end{array}
\]
It follows that the double extension of $\fbrj (2;3)$ is a~\textit{pre-Lie superalgebra}, see Appendix.
\sssbegin{Claim} For $\fa=\mathfrak{sl}(1|2)$, $\mathfrak{osp}(2|3)$, $\mathfrak{osp}(1|4)$ for $p>2$, the Brown algebra $\mathfrak{br}(2,\varepsilon)$ for $p=3$ and $\fbrj (2;5)$ for $p=5$, we have $\mathrm{H}^1 (\fa; \fa)=~0$.
\end{Claim}
Therefore, these Lie superalgebras do not have non-trivial double extentions. 
\subsection{Manin triples related with $\fhei(2)$ for $p=2$} \label{MT}
\subsubsection{$2$-structures on Manin triples.}\label{Maninp=2}
Let $(\fh, [2]_\fh)$ be a~finite-dimensional restricted Lie algebra (not necessarily \lq\lq NIS\rq\rq), and let $\fh^*$ have the structure of an abelian Lie algebra. Set $\fa:=\fh\oplus \fh^*$, and define the bracket of two elements on $\fa$ as follows:
\begin{equation}
\label{bracketh}
[h+\pi, h'+\pi']_{\fa}:=[h,h']_\fh+\pi\circ \text{ad}_{h'}+\pi'\circ \text{ad}_{h}\quad \text{for any $h+\pi,  \ h'+\pi'\in\fa$.}
\end{equation}
It is easy to show that the bracket $[\cdot , \cdot]_{\fa}$ defined by Eq. (\ref{bracketh}) satisfies the Jacobi identity.

Define the $2$-structure on $\fa$ as follows (for any $h\in \fh$ and $\pi\in \fh^*$, hence for any $h+\pi \in \fa$):
\begin{equation}
\label{squaringh}
(h+\pi)^{[2]_\fa}:=h^{[2]_\fh}+\pi\circ \text{ad}_h.
\end{equation}
Let us show that the mapping defined by Eq. (\ref{squaringh}) is a~$2$-mapping. 
 
 Indeed, 
\[
\begin{array}{lcl}
[h+\pi,[h+\pi, h'+\pi']_{\fa}]_\fa&=&[ h+\pi,[h,h']_\fh+\pi\circ \text{ad}_{h'}+\pi'\circ \text{ad}_{h}]_\fa\\[2mm]
&=&[ h,[h,h']_\fh]_\fh+\pi\circ \ad_{[h,h']_\fh} + (\pi\circ \text{ad}_{h'}+\pi'\circ \text{ad}_{h})\circ \ad_{h}\\[2mm]
&=&[h^{[2]_\fh} ,  h']_\fa+\pi\circ \ad_{[h,h']_\fh}+\pi\circ \text{ad}_{h'}\circ \ad_{h}+\pi'\circ \text{ad}_{h}\circ \ad_{h}\\[2mm]
&=&[h^{[2]_\fh} ,  h']_\fa+\pi\circ \text{ad}_{h}\circ \ad_{h'}+\pi'\circ \text{ad}_{h^{[2]_\fh}}\\[2mm]
&=&[(h+\pi)^{[2]}, h'+\pi']_\fa.
\end{array}
\]
On the other hand,
\[
\begin{array}{lcl}
(h+\pi+h'+\pi')^{[2]}&=&(h+h')^{[2]_\fh}+(\pi+\pi')\circ \text{ad}_{h+h'}\\[2mm]
&=&h^{[2]_\fh}+h'^{[2]_\fh}+[h,h']_\fh+\pi\circ \text{ad}_{h}+\pi\circ \text{ad}_{h'}
+\pi'\circ \text{ad}_{h}
+\pi'\circ \text{ad}_{h'}\\[2mm]
&=&(h+\pi)^{[2]}+(h'+\pi')^{[2]}+[h+\pi,h'+\pi']_\fg.
\end{array}
\]
Define a~bilinear form on $\fa$ as follows:
\begin{equation}\label{Mform}
\mathscr{B}_{\fa}(h+\pi, h'+\pi'):=\pi(h')+\pi'(h) \quad \text{for any $h+\pi, \ h'+\pi'\in\fa$.}
\end{equation}
It is easy to show that the bilinear form $\mathscr{B}_\fa$ is NIS. 
\subsubsection{Heisenberg algebra $\fhei(2)$}
Consider the Heisenberg algebra $\fhei(2)$ spanned by $p,q$ and $z$, with the only nonzero bracket: $[p,q]=z$. Let us consider a~2-structure given by
\[
p^{[2]}=q^{[2]}=0, \quad z^{[2]}=z.
\]
We consider the NIS-algebra $\fa:=\fhei(2)\oplus \fhei(2)^*$ constructed as in \S \ref{Maninp=2}.  A direct computation using Eq. (\ref{squaringh}) shows that (for any $s,w,u,v\in \mathbb{K}$)
\begin{eqnarray*} 
(rz+s p+wq+up^*+v q^*+t z^*)^{[2]_\fa}=(r^2+sw)z+st q^*+wt p^*.
\end{eqnarray*}
A direct computation using Eqs. (\ref{squaringh}) and (\ref{bracketh}) shows that the only nonzero brackets are 
\[
[p,q]_\fa = z,\quad [p,z^*]_\fa = q^*,\quad [q, z^*]_\fa=p^*.
\]

\sssbegin{Claim}
\label{der1}
The space $\mathrm{H}^1_{\mathrm{res}}(\fa;\fa)$ is spanned by the (classes of the) following cocycles:
\[
\begin{array}{lcllcllcllcl}
\mathscr{D}_1&=&q^*\otimes  \widehat{p}, &\mathscr{D}_2&=&q^*\otimes  \widehat{q}, &\mathscr{D}_3&=&q^*\otimes  \widehat{z^*},\\[2mm] 
\mathscr{D}_4&=&p^*\otimes  \widehat{p},&\mathscr{D}_5&=&p^*\otimes  \widehat{z^*}, &\mathscr{D}_6&=&z\otimes  \widehat{z^*},\\[2mm]
\mathscr{D}_7&=&p\otimes  \widehat{p}+q^*\otimes  \widehat{q^*}+z\otimes  \widehat{z},& \mathscr{D}_8&=&q\otimes  \widehat{q}+p^*\otimes  \widehat{p^*}+z\otimes  \widehat{z},&&& \\[2mm]
 \mathscr{D}_{9}&=&q^*\otimes  \widehat{q^*}+p^*\otimes  \widehat{p^*}+z^*\otimes  \widehat{z^*}.
\end{array}
\]
\end{Claim}

Let us fix an ordered basis as follows: $p,q,z, p^*, q^*, z^*$. In this basis, the Gram matrix of the bilinear form $\mathscr{B}_\fa$ in (\ref{Mform}) is 
\[
\mathrm{antidiag}(I_3, I_3), \text{ where $I_n$ denotes the  $n\times n$ unit matrix.}
\]
\sssbegin{Proposition}
Up to an isomorphism, $\fa$ admits a~unique non-trivial double extension given by the derivation $\mathscr{D}_7+\mathscr{D}_9.$
\end{Proposition}
\begin{proof}
Any derivation $\mathscr{D}$ has the following representation as a~supermatrix in the standard format:
\[
\mathscr{D}=\left (
\begin{array}{c|c}
A&B\\
\hline
C& F
\end{array}
\right ).
\] It follows that $\mathscr{B}_\fa$ is $\mathscr{D}$-invariant if and only if $F=A^t, B^t=B$ and $C^t=C$. Let us consider the most general derivation $\mathscr{D}=\mathop{\sum}\limits_{1 \leq i\leq 9}\mu_i \mathscr{D}_i$, where $\mathscr{D}_i$ are the cocycles given in Claim~\ref{der1}. In the same basis $p,q,z,p^*,q^*,z^*$, we have
\small{
\[
\mathscr{D}=\left (
\begin{array}{c|c}
\mu_7 E^{1,1}+(\mu_7+\mu_8) E^{3,3}+\mu_8 E^{2,2}&\mu_6 E^{3,3}\\
\hline
\mu_4 E^{1,1}+\mu_1 E^{2,1}+\mu_2 E^{2,2}&\mu_8 E^{1,1}+\mu_7 E^{2,2}+\mu_3 E^{2,3}+\mu_5 E^{1,3}+\mu_9 I
\end{array}
\right ),
\]
}
where $E^{i,j}$ is the $(i,j)$th $3\times 3$ matrix unit. 

It follows that $\mathscr{B}_\fa$ is $\mathscr{D}$-invariant if and only if $\mu_1=\mu_3=\mu_5=0$ and $\mu_9=\mu_7+\mu_8$. 

The most general derivation is of the form 
\[
\mu_2\mathscr{D}_2+\mu_4 \mathscr{D}_4+\mu_6 \mathscr{D}_6+\mu_7 \mathscr{D}_7+\mu_8\mathscr{D}_8+(\mu_{7}+\mu_8)\mathscr{D}_{9}.
\] 
Since $p=2$, we also need to check that $\mathscr{B}_\fa(a,\mathscr{D}(a))=0$ for all $a\in \fa$. Let 
\[
a=rz+s p+wq+up^*+v q^*+t z^*\in \fg.
\] 
The fact that
\[
\begin{array}{lcl} 
\mathscr{B}_\fa(a,\mathscr{D}(a))&=&\mathscr{B}_\fa(a,(\mu_2 w+\mu_8 v)q^*+(\mu_4s+u\mu_7)p^*+(\mu_6t+\mu_7r+\mu_8r)z\\[2mm]
&&+\mu_7 sp+\mu_8wq+(\mu_7+\mu_8)tz^*)\\[2mm]
&=&(\mu_2 w)w+(\mu_4s)s+(\mu_6t)t=0,
\end{array}
\] 
implies that $\mu_2=\mu_4=\mu_6=0$. 

Now, we define two quadratic forms on $\fa$ as follows: 
\[
\begin{array}{clc}
\mathscr{P}_1(rz+s p+wq+t z^*+up^*+v q^*)&=&rt+su,\\[2mm]
\mathscr{P}_2(rz+s p+wq+t z^*+up^*+v q^*)&=&rt+wv.
\end{array}
\] 
Now, let 
\[
b=rz+s p+wq+up^*+v q^*+t z^*,\; a=\tilde rz+\tilde s p+\tilde wq+\tilde up^*+\tilde v q^*+\tilde t z^* \in \fa.
\]

Direct computations show that 

\[
\begin{array}{lcl}
\mathscr{B}_\fa(a, \mathscr{D}(b))&=&\mathscr{B}_\fa(a, \mu_8v q^*+ u\mu_7 p^*+(\mu_7r+\mu_8r)z+\mu_7 sp+\mu_8wq+(\mu_7+\mu_8)tz^*)\\[2mm]
&=&(\mu_8v)\tilde w+(u\mu_7)\tilde s+(\mu_7r+\mu_8r)\tilde t+\mu_7 s\tilde u+\mu_8w\tilde v+(\mu_7+\mu_8)t\tilde r\\[2mm]
&=&\mu_7 (s\tilde u+u\tilde s)+(\mu_7+\mu_8)(r \tilde t+t\tilde r)+\mu_8(w\tilde v+v\tilde w).
\end{array}
\]

Hence, we have two $\mathscr{D}$-extensions given by the following data (where $\alpha, \beta \in \mathbb{K}$):
\begin{equation}\label{g1a0}
\begin{array}{l}
(\mathscr{D}=\mathscr{D}_7+\mathscr{D}_9,\ \  a_0=\beta q^*,\ \   b_0=\alpha q^*, \ \  \mathscr{P}_1, \ \  m,\ \  l,\ \  \gamma=1),\\[3mm] 
(\mathscr{D}=\mathscr{D}_8+\mathscr{D}_9, \ \  a_0=\beta p^*, \ \  b_0=\alpha p^*,\ \  \mathscr{P}_2, \ \  \tilde m, \ \  \tilde l,\ \   \tilde \gamma=1 ).
\end{array}
\end{equation}

Let us show that these two $\mathscr{D}$-extensions are isomorphic if $m,l,\tilde m$ and $\tilde l$ are suitably chosen. Indeed, the isomorphism is given by (for notation, see Theorem \ref{Isom1})
\[
\begin{array}{rclrclrclrcl}
\pi_0(z)&=& z, & \pi_0(z^*)&=&z^*, &  \pi_0(p)&=&q, \\[2mm]
 \pi_0(q)&=&p & \pi_0(q^*)&=&p^*,& \pi_0(p^*)&=&q^*,\\[2mm]
 \varkappa&=&0, & \lambda&=&1,& \rho&=&0.
\end{array}
\]

On the other hand, let us show that the $\mathscr{D}$-extension by means of $\mathscr{D}_7+\mathscr{D}_9$ is not a~trivial one; namely, it is not isomorphic to the one by means of $\ad_T$ for some $T\in \fa$. Suppose there is such an isomorphism, say $\pi$. Let us write 
\[
\pi_0(z)=mz, \quad \pi_0(z^*)=m_1z+m_2p+m_3q+m_4p^*+m_5q^*+m^{-1}z^*.
\]
Now, because $q^*=[z^*,p]$, it follows that for some $c_1, c_2, c_3 \in\mathbb{K}$, we have
\[\pi_0(q^*)=[m_1z+m_2p+m_3q+m_4p^*+m_5q^*+m^{-1}z^*, \pi_0(p)]=c_1p^*+c_2q^*+c_3z.
\] 
Similarly, since $p^*=[z^*,q]$, it follows that for some $\tilde c_1, \tilde c_2, \tilde c_3\in \mathbb{K}$, we have 
\[
\pi_0(p^*)=[m^{-1}z^*, \pi_0(q)]=\tilde c_1p^*+\tilde c_2q^*+\tilde c_3 z.
\] 
Let $T=W_1p+W_2p^*+W_3q+W_4q^*+W_5z+W_6z^*$, where $W_i \in \mathbb{K}$. We have
\[
\begin{array}{ll}
&((\mathscr{D}_7+\mathscr{D}_9)\circ \pi_0-\pi_0\circ \ad_T)(z^*)\\
=&(\mathscr{D}_7+\mathscr{D}_9)(m_1z+m_2p+m_3q+m_4p^*+m_5q^*+m^{-1}z^*)
-\pi_0 [T, z^*]\\[2mm]
=&m^{-1}z^*+m_1z+m_2p+m_4p^*- W_1\pi_0(q^*)-W_3\pi_0(p^*)\\[2mm]
=&m^{-1}z^*+m_1z+m_2p+m_4p^*- W_1(c_1p^*+c_2q^*+c_3 z)\\[2mm]
&-W_3 (\tilde c_1p^*+\tilde c_2q^*+\tilde c_3 z).
\end{array}
\] 
But this is never zero, hence a~contradiction.
\end{proof}
\subsection{Vectorial Lie (super)algebras}\label{vectorial}
Over any field $\mathbb{K}$ of characteristic $p > 0$, consider not polynomial coefficients but divided powers in $n$ indeterminates, whose powers are bounded by the shearing vector $\underline{N}= (N_1, ...,N_n)$. We get a~commutative algebra (here $p^\infty := \infty)$
\[
\mathscr{O}(n;\underline{N}):=\mathbb{K} \left [u;\underline{N}\right ]:=\mathrm{Span}_{\mathbb{K}} ( u^{(r)} \,|\, 0 \leq r_i < p^{N_i}) ,
\]
where $u^{(r)} = \prod_{1\leq i \leq n} u^{(r_i)}_i .$ The addition in $\mathscr{O}(n;\underline{N})$ is natural; the multiplication is defined by 
\[
u^{(r_i)}_i \cdot u^{(s_i)}_i= \binom{r_i+s_i}{r_i} u^{(r_i+s_i)}_i. 
\]
Set ${\bf 1} := (1, . . . ,1)$ and set $\tau(\underline{N}):=(p^{N_1}-1,\ldots, p^{N_n}-1)$ often abbreviated to $\tau$.

The \textit{distinguished} partial derivatives $\partial_i$, each of them serving as several partial
derivatives at once, for each of the generators $u_i , u^{(p)}_i , u^{(p^2)}_i , \cdots$ (or, in terms of $y_{i,j} :=u_i^{(p^j-1)})$ are defined by the formula
\[
\partial_i(u^{(k)}_j ) := \delta_{ij}u^{(k-1)}_j \text{ for all $k$, i.e., $\partial_i =\sum_{j\geq 1} (-1)^{j-1}y^{p-1}_{i,1} \cdots y^{p-1}_{i,j-1}\partial _{y_{i,j}}.$}
\]


\subsubsection{$\mathfrak{vect}(n;\underline{N})$} 

The \textit{general vectorial Lie algebra}, known as the \textit{Jacobson-Witt algebra}:
\[
\begin{array}{lcl}
\mathfrak{vect}(n;\underline{N})&:=&\left \{\sum_i f_i \partial_i \;|\; f_i \in \mathscr{O} (n;\underline{N}) \right \},
\end{array}
\]
where the bracket is given by the Lie bracket of vector fields. The Lie algebra $\mathfrak{vect}(n;\underline{N})$ has a~NIS (\cite[Theorems 6.3 and 6.4 in Ch. 4]{SF}) if and only if either $n = 1$ and $p = 3$, in which case NIS is
\begin{equation}\label{vect1p3}
\mathscr{B}(u^{(a)}\partial, u^{(b)}\partial) := \int 
u^{(a)}u^{(b)}du,
\end{equation}
or $n = p = 2$, in which case NIS is
\[
\mathscr{B}(u^{(a)}\partial_i,u^{(b)}\partial_j ):=(i+j)\int u^{(a)}u^{(b)} du_1\wedge du_2,
\]
where $\int f (u)du_1\wedge \cdots \wedge du_n := \text{coefficient of } u^{(\tau (\underline{N}))}$ in the Taylor series expansion of $f(u)$.

\parbegin{Remark} It was shown in \cite{BKLS} that there is no NIS on the simple Lie
algebras denoted there by $\mathfrak{svect}_{exp_i} (n;\underline{N})$ and $\mathfrak{svect}^{(1)}
_{1+\bar{u}}(n;\underline{N})$ --- deforms of the divergence-free subalgebra of $\fvect(n;\underline{N})$. 
For examples of NISes on simple $\mathbb{Z}$-graded vectorial Lie algebras in characteristic $p > 0$, see \cite{Dz, Fa}; these examples are reproduced in \cite{BKLS}.
\end{Remark}

Since we are interested only in restricted Lie (super)algebras which can only exist for $\underline{N}={\bf 1}$, we mostly skip $\underline{N}$.

\parbegin{Claim} \textup{(}i\textup{)}   \cite[Theorem 8.5]{SF} For $p=3$, $ \mathrm{H}^1(\mathfrak{vect}(1);\mathfrak{vect}(1))=0$.  Hence, no non-trivial double extensions of $\mathfrak{vect}(1)$. 

\textup{(}ii\textup{)}  For $p=2$, $ \mathrm{H}^1(\mathfrak{vect}(2);\mathfrak{vect}(2))=0$. Hence, no  non-trivial double extensions of $\mathfrak{vect}(2)$. 
\end{Claim}

\subsubsection{$\mathfrak{svect}(3;\underline{N})$ for $p=2$} Set
\[
\begin{array}{lcl}
\mathfrak{svect}(n;\underline{N})&:=&\left \{\mathscr{D}  \in \mathfrak{vect}(n;\underline{N})\;|\; \mathrm{div}(\mathscr{D})=0 \right \}.
\end{array}
\]
This Lie algebra is not simple; however its first derived subalgebra $\mathfrak{svect}^{(1)}(n;\underline{N})$ is simple for $n\geq 3$, see \cite{S, SF}. Consider the maps 
\[
D_{i,j}: \mathscr{O}(n;\underline{N}) \rightarrow \mathfrak{vect}(n;\underline{N}), \quad f\mapsto D_{i,j}(f)=\partial_j(f)\partial_i-\partial_i(f)\partial_j.
\]
Clearly, 
\[
\mathfrak{svect}^{(1)}(n;\underline{N})=\Span \{D_{i,j}(f)\; | \; f\in \mathscr{O}(n,\underline{N}), 1\leq i<j\leq n \}.
\]
The Lie algebra $\mathfrak{svect}^{(1)}(n;\underline{N})$ has a~NIS if and only if $n = 3$; explicitly, we set 
\begin{equation}\label{Gramsvect}
\mathscr{B} (\partial_i, D_{j,k}(u^{(\tau (\underline{N}))})) = \text{sign}(i, j, k), 
\end{equation}
and extend the form $\mathscr{B}$ to other pairs of elements by invariance and linearity.

We have the following exceptional isomorphisms 
\[
\mathfrak{svect}^{(1)}(3; {\bf 1})\simeq \mathfrak{psl}(4)\simeq \mathfrak{h}_\omega^{(1)}(4;{\bf 1}) \quad \text{for $p = 2$ (as shown in \cite{CK})}.
\] 
The case of the Hamiltonian Lie superalgebra $\fh_\omega$ preserving the symplectic form $\omega$ with constant coefficients is completely investigated in \cite{BLS} for all possible types of $\omega$, see Subsection~\ref{Secpsl4}.

\parbegin{Claim} $\dim \mathrm{H}^1_{\mathrm{res}}(\mathfrak{svect}^{(1)}(3);\mathfrak{svect}^{(1)}(3))=2$ is spanned by the cocycles $\mathscr{D}_3$ of degree $3$ and $\mathscr{D}_0$ of degree $0$. 
\tiny
\[
\begin{array}{ccl}
\deg=3: \mathscr{D}_3&=& 2 D_{1,3}(u^{(\tau)}) \otimes \widehat{(D_{2,3}(u^{(\tau-\epsilon_2-2\epsilon_3)}) )}+
D_{1,2}(u^{(\tau)}) \otimes \widehat{(D_{2,3}(u^{(\tau-\epsilon_3-2\epsilon_2)}) )}\\[2mm]
&&+
\left ( 2D_{1,3}(u^{(\tau-\epsilon_1)}) +2D_{1,2}(u^{(\tau-\epsilon_1)}) \right )\otimes \widehat{(D_{1,3}(u^{(\tau-2\epsilon_2-2\epsilon_3)}) )},
\\[2mm]
&&
+D_{1,3}(u^{(\tau-2\epsilon_1)}) \otimes \widehat{\partial_3}
+D_{1,2}(u^{(\tau-2\epsilon_1)}) \otimes \widehat{\partial_2}\\[2mm]
\deg=0:\mathscr{D}_0&=& 
\left ( D_{1,2}(u^{(\tau-\epsilon_1)})+2D_{1,3}(u^{(\tau-\epsilon_1)})+2D_{2,3}(u^{(\tau-\epsilon_1)}) \right )\otimes \widehat{(D_{1,2}(u^{(\tau-\epsilon_1)}) )}\\[2mm]
&&
+\left ( 2D_{1,3}(u^{(\tau-\epsilon_2)}) +D_{1,2}(u^{(\tau-\epsilon_2)}) \right ) \otimes \widehat{(D_{1,2}(u^{(\tau-\epsilon_2)}) )}
+2D_{1,3}(u^{(\tau)}) \otimes \widehat{(D_{1,2}(u^{(\tau)}) )} \\[2mm]
&&
+D_{2,3}(u^{(\tau-2\epsilon_1)}) \otimes \widehat{(D_{1,2}(u^{(\tau-2\epsilon_1)}) )} 
+D_{1,3}(u^{(\tau-\epsilon_1-\epsilon_2)}) \otimes \widehat{(D_{1,2}(u^{(\tau-\epsilon_1-\epsilon_2)}) )}\\[2mm]
&&
+2D_{2,3}(u^{(\tau-\epsilon_2)}) \otimes \widehat{(D_{2,3}(u^{(\tau-\epsilon_2)}) )}

+2D_{1,3}(u^{(\tau-\epsilon_2-\epsilon_3)}) \otimes \widehat{(D_{1,3}(u^{(\tau-\epsilon_2-\epsilon_3)}) )}\\[2mm]
&&
+2D_{2,3}(u^{(\tau-\epsilon_2-\epsilon_3)}) \otimes \widehat{(D_{2,3}(u^{(\tau-\epsilon_2-\epsilon_3)}) )}
+D_{1,2}(u^{(\tau-\epsilon_2-\epsilon_3)}) \otimes \widehat{(D_{1,2}(u^{(\tau-\epsilon_2-\epsilon_3)}) )}\\[2mm]
&&
+\left ( D_{1,3}(u^{(\tau-2\epsilon_2)}) +D_{2,3}(u^{(\tau-2\epsilon_2)}) \right )\otimes \widehat{(D_{1,3}(u^{(\tau-2\epsilon_2)}) )}
+2\partial_1\otimes \partial_1+ 2\partial_1\otimes \partial_2\\[2mm]
&&
+D_{1,3}(u^{(\tau-2\epsilon_1)}) \otimes \widehat{(D_{1,3}(u^{(\tau-2\epsilon_1)}) )}
+2D_{1,3}(u^{(\tau-\epsilon_1-\epsilon_3)}) \otimes \widehat{(D_{1,3}(u^{(\tau-\epsilon_1-\epsilon_3)}) )}\\[2mm]
&&
+2D_{2,3}(u^{(\tau-\epsilon_1-\epsilon_3)}) \otimes \widehat{(D_{2,3}(u^{(\tau-\epsilon_1-\epsilon_3)}) )}
+D_{1,2}(u^{(\tau-\epsilon_2-2\epsilon_3)}) \otimes \widehat{(D_{1,2}(u^{(\tau-\epsilon_2-2\epsilon_3)}) )}\\[2mm]
&&
+D_{1,3}(u^{(\tau-2\epsilon_2-\epsilon_3)}) \otimes \widehat{(D_{1,3}(u^{(\tau-2\epsilon_2-\epsilon_3)}) )}
+D_{2,3}(u^{(\tau-2\epsilon_1-\epsilon_3)}) \otimes \widehat{(D_{2,3}(u^{(\tau-2\epsilon_1-\epsilon_3)}) )}\\[2mm]
&&
+D_{1,3}(u^{(\tau-\epsilon_1\epsilon_2-\epsilon_3)}) \otimes \widehat{(D_{1,3}(u^{(\tau-\epsilon_1-\epsilon_2-\epsilon_3)}) )}
+2D_{2,3}(u^{(\tau-\epsilon_2-2\epsilon_3)}) \otimes \widehat{(D_{2,3}(u^{(\tau-\epsilon_2-2\epsilon_3)}) )}\\[2mm]
&&
+2D_{1,3}(u^{(\tau-\epsilon_1-2\epsilon_3)}) \otimes \widehat{(D_{1,3}(u^{(\tau-\epsilon_1-2\epsilon_3)}) )}
+2D_{1,2}(u^{(\tau-\epsilon_1-2\epsilon_2)}) \otimes \widehat{(D_{1,2}(u^{(\tau-\epsilon_1-2\epsilon_2)}) )}\\[2mm]
&&
+2D_{2,3}(u^{(\tau-\epsilon_1-2\epsilon_3)}) \otimes \widehat{(D_{2,3}(u^{(\tau-\epsilon_1-2\epsilon_3)}) )}
+D_{2,3}(u^{(\tau-2\epsilon_1-2\epsilon_3)}) \otimes \widehat{(D_{2,3}(u^{(\tau-2\epsilon_1-2\epsilon_3)}) )}\\[2mm]
&&
+D_{2,3}(u^{(\tau-2\epsilon_2-\epsilon_3)}) \otimes \widehat{(D_{2,3}(u^{(\tau-2\epsilon_2-\epsilon_3)}) )}
+D_{1,3}(u^{(\tau-2\epsilon_1-\epsilon_3)}) \otimes \widehat{(D_{1,3}(u^{(\tau-2\epsilon_1-\epsilon_3)}) )}\\[2mm]
&&
+2D_{1,2}(u^{(\tau-\epsilon_2-2\epsilon_1)}) \otimes \widehat{(D_{1,2}(u^{(\tau-\epsilon_2-2\epsilon_1)}) )}+D_{1,3}(u^{(\tau-2\epsilon_2-2\epsilon_3)}) \otimes \widehat{(D_{1,3}(u^{(\tau-2\epsilon_2-2\epsilon_3)}) )}\\[2mm]
&&
+D_{1,3}(u^{(\tau-2\epsilon_2-2\epsilon_1)}) \otimes \widehat{(D_{1,3}(u^{(\tau-2\epsilon_2-2\epsilon_1)}) )}
+D_{2,3}(u^{(\tau-2\epsilon_2-2\epsilon_1)}) \otimes \widehat{(D_{2,3}(u^{(\tau-2\epsilon_2-2\epsilon_1)}) )}.
\end{array}
   \]
   \end{Claim}
   
 A direct computation shows that the bilinear form given in Eq. (\ref{Gramsvect}) is {\bf not} $\mathscr{D}_3$-invariant. Therefore, the Lie algebra $\mathfrak{svect}^{(1)}(3)$ cannot be double extended by means of $\mathscr{D}_3$. However, a~direct computation shows that this bilinear form is $\mathscr{D}_0$-invariant and, moreover, $(\mathscr{D}_0)^3~=~0$ (cube of the operator). Therefore, $a_0=\gamma=0$, see \S \ref{Rod}.  
 
 Now we define the cubic form $ \mathscr{P}(\sum_{1\leq i \leq 52}\lambda_ie_i)$ to be equal to
 \begin{equation}
 \label{cubicsvect}
 \begin{array}{l}
2 \lambda _{32} \lambda _2^2+\lambda _9 \lambda _{12}
   \lambda _2+\lambda _{10} \lambda _{17} \lambda _2 +2
   \lambda _4 \lambda _{20} \lambda _2+2 \lambda _8
   \lambda _{21} \lambda _2\\[2mm]
   +\lambda _3 \lambda _{37}
   \lambda _2 +\lambda _7 \lambda _8^2  +\lambda _4^2 \lambda _9 +\lambda _4 \lambda _8 \lambda _{10} 
 +2
   \lambda _4 \lambda _8 \lambda _{11}
   +2 \lambda _1
   \lambda _8 \lambda _{12}\\[2mm]
   +\lambda _3 \lambda _{11}
   \lambda _{12} +\lambda _3 \lambda _8 \lambda_{15}+\lambda _1 \lambda _4 \lambda _{17} +2 \lambda
   _3 \lambda _7 \lambda _{17}
   +2 \lambda _3 \lambda _4
   \lambda _{22}+2 \lambda _3^2 \lambda _{28},
   \end{array}
 \end{equation}
where it suffices to describe the $e_i$ that appear in the expression of $\mathscr{P}(a)$:
\[
\begin{array}{lll}
e_2=\partial_2, & e_3=\partial_3, &  e_4=2D_{1,2}(u^{(\tau-2\epsilon_2-2\epsilon_3)}),\\[2mm]
e_7=D_{2,3}(u^{(\tau-2\epsilon_1-2\epsilon_2)}), & e_8=D_{1,3}(u^{(\tau-2\epsilon_2-2\epsilon_3)}), & e_9=D_{2,3}(u^{(\tau-2\epsilon_1-2\epsilon_3)}),\\[2mm] 
e_{10}=D_{1,2}(u^{(\tau-\epsilon_1-\epsilon_2-2\epsilon_3)}),&
e_{11}=D_{1,3}(u^{(\tau-2\epsilon_2-\epsilon_3-\epsilon_1)}),& 
e_{12}=D_{2,3}(u^{(\tau-\epsilon_3-2\epsilon_2)}), \\[2mm]
e_{15}=2D_{1,2}(u^{(\tau-\epsilon_3-2\epsilon_2)}),&
e_{17}=D_{2,3}(u^{(\tau-2\epsilon_3-\epsilon_2)}),&
e_{20}=D_{1,3}(u^{(\tau-2\epsilon_3-\epsilon_2)}),\\[2mm]
e_{21}=D_{1,2}(u^{(\tau-2\epsilon_3-\epsilon_2)}),&
e_{22}=D_{1,3}(u^{(\tau-\epsilon_3-2\epsilon_2)})&
e_{28}=D_{2,3}(u^{(\tau-2\epsilon_2)}),\\[2mm]
e_{32}=D_{2,3}(u^{(\tau-2\epsilon_3)}),&
e_{37}=D_{2,3}(u^{(\tau-\epsilon_3-\epsilon_2)}).
\end{array}
\]

Let us summarize and introduce notation:
\begin{equation}\label{svect}
\footnotesize
\renewcommand{\arraystretch}{1.4}
\begin{tabular}{|c|c|c|c|c|c|} \hline
Derivation & $\mathscr{D}$-invariance &$\mathscr{P}(a)$&$\gamma$ &$a_0$&
Double extension\\
\hline
$\mathscr{D}_3$ & No &---&---&---& --- \\ \hline
$\mathscr{D}_0$ & Yes & see Eq. (\ref{cubicsvect}) &0& 0&$\widehat \fsvect^{(1)}(3)$\\ \hline
\end{tabular}
\end{equation}

\subsubsection{The Hamiltonian Lie (super)algebras $\fh^{(1)}_\omega(a;{\bf 1}|b)$ for $a+b>3$ and $p=2$} \label{Secpsl4} 

There are four types of symplectic forms $\omega$ with constant coefficients, 
as  Lebedev showed, see \cite{LeD}. For the double extensions of $\fh^{(1)}_\omega(a;{\bf 1}|b)$ for $a+b>3$ and $p=2$, see \cite{BLS}. 

For $p>2$ and in the non-super case, Skryabin classified the normal shapes of symplectic forms $\omega$ with non-constant coefficients, see \cite{Sk, Sk1}; some of these shapes exist for $p=2$. If $p=2$, no classification is known; for new examples, see \cite{KoKCh}. Some of these $\mathfrak{h}_\omega(2n;\underline{N})$ have a~NIS, see \cite{BKLS}.  To consider these cases, as well as deformations of other simple vectorial Lie algebras and superalgebras is an open problem.

\subsubsection{The exceptional Lie superalgebra  $\fvas$ for $p=3$, see \cite{BKLS}, has no double extensions. Sketch of the proof}\label{vas} The Lie superalgebra $\fg:=\fvas$ has a grading operator $D$, hence the central extension should be graded. If, moreover, we wish this central extension be inherited from $\fsle^{(1)}(4)$, see \cite{Sh14}, the central element must be odd of degree $-2$. This means that the bracket should define projectively $\fg_0$-invariant 2-form $\omega$ on $\fg_{-1}$.

If $\omega$ were degenerate, its kernel would have been $\fg_0$-invariant subspace of $\fg_{-1}$, but the latter is irreducible. 

Therefore, $\omega$ is a non-degenerate odd projectively $\fg_0$-invariant bilinear form on $\fg_{-1}$. Hence, $\fg_0\subset \fpe(4)\ltimes\mathbb{K} D$. Precisely as 0th component of $\fle(4)$, see \cite{Le}, but $\fas=\fvas_0$ (for its definition, see \cite{Sh14}) can not be embedded into $\fpe(4)\ltimes\mathbb{K} D$.

\subsection{$\fosp(1|2)$, $\mathfrak{g}^{(1)}(2,3)/\fz$ and $\mathfrak{g}^{(1)}(3,3)/\fz$ for $p=3$}
\subsubsection{$\fosp(1|2)$} \label{preosp}
Let us fix a~basis of $\mathfrak{osp}(1|2)$ generated by the root vectors $x_1, x_2=[x_1, x_1]$ (positive), $y_1, y_2=[y_1, y_1]$ (negative), and $h_1=[x_1,y_1]$. The Lie superalgebra $\mathfrak{osp}(1|2)$ admits a~NIS given in the ordered basis $e_1=[x_1, y_1]$, $e_2= x_1$, $e_3=x_2$, $e_4=y_1$, $e_5=y_2$ by the Gram matrix (here $E^{i,j}$ is the $(ij)$-th matrix unit)
\[
\mathscr{B}_{\mathfrak{osp}(1|2)}=2 E^{1,1}+2 E^{2,4}+2E^{3,5}+E^{4,2}+2 E^{5,3}.
\]
There is a~$3|6$-structure on $\mathfrak{osp}(1|2)$ given by 
\[
e_1^{[3]}=e_1, \quad  e_2^{[6]}=e_4^{[6]}=e_3^{[3]}=e_5^{[3]}=0.
\]

\parbegin{Claim} \textup{(\cite{BGLL1})} \textup{(i)} For $p=3$, the  space $\mathrm{H}^1 (\fosp(1|2); \fosp(1|2))$ is spanned by the odd cocycles: 
\begin{equation}\label{cocyposp}\tiny
\begin{array}{ll}
\deg=-3:&\mathscr{D}_{-3}=2 y_1 \otimes  \widehat{x_2}+
   y_2 \otimes \widehat{x_1}, \\[2mm]

\deg=3: &\mathscr{D}_{3}=x_1\otimes  \widehat{y_2}+
   x_2\otimes  \widehat{y_1}.
\end{array}
\end{equation}

\textup{(ii)}  For $p>3$,  $\mathrm{H}^1 (\fosp(1|2); \fosp(1|2))=0$.
\end{Claim}
A long, but easy, computation shows that $\mathscr{B}_{\fosp (1|2)}$ is both $\mathscr{D}_{-3}$-invariant and $\mathscr{D}_{3}$-invariant. Moreover, it is easy to show that $(\mathscr{D}_{-3})^2=(\mathscr{D}_{3})^2=0$ (squares of the operators). Therefore, conditions (\ref{DB0}) are satisfied because $\fosp(1|2)$ is simple. However, the condition (\ref{p=3Con1}) is violated. Indeed, for the derivation $\mathscr{D}_{-3}$ we have
\[
\mathscr{B}_{\mathfrak{osp}(1|2)}(\mathscr{D}_{-3}(x_1),[x_1,x_1])=\mathscr{B}_{\mathfrak{osp}(1|2)}(2 y_2, x_2)=1\not=0.
\]
Therefore, the double extension of $\mathfrak{osp}(1|2)$ is a~pre-Lie superalgebra, see Appendix. The isomorphism $\pi$ defined by   
\[
\pi(x_1)= 2y_1, \; \pi(y_1)=x_1,
\]
shows that the derivations $\mathscr{D}_{-3}$ and $\mathscr{D}_{3}$ produce the same pre-Lie superalgebra.
\subsubsection{$\fg(2,3)^{(1)}/\fz$, see \cite{BGL}}
Consider the Lie superalgebra  $\fg(2,3)$ with the Cartan matrix
\begin{equation*}\label{eq2.3.2b}\tiny
\begin{pmatrix}
 2 & -1 & -1 \\
 -1 & 2 & -1 \\
 -1 & -1 & 0\end{pmatrix} \text{ with positive roots }{\footnotesize
 \begin{array}{l}
x_1,\;\;x_2,\;\;x_3,\\
x_4=[x_1, x_2],\;\;
x_5=[x_1, x_3],\;\;x_6=[x_2, x_3],\\
x_7=[x_3, [x_1, x_2]],\\
x_8=[[x_1, x_2], [x_1, x_3]],\;\; x_9=[[x_1,
x_2], [x_2, x_3]],\\
x_{10}=[[x_1, x_2], [x_3, [x_1, x_2]]],\\
x_{11}=[[x_3, [x_1, x_2]], [x_3, [x_1, x_2]]]\end{array}
}
\end{equation*} 
Denote by $y$'s the corresponding negative roots, and set $h_i = [x_i,y_i]$ for $i=1,2,3$.

The simple core $\fa:=\fg(2,3)^{(1)}/\fz $ admits a~NIS given in the ordered basis 
\begin{equation}
\label{basg23}
h_2:=[x_2,y_2], \; h_3:= [x_3,y_3],\;  x_1,\ldots x_{11},\; y_1, \ldots,y_{11},
\end{equation}
by the Gram matrix
\[
\begin{array}{l}
\mathscr{B}_{\fa}=\left (
\begin{array}{ccc}
S& 0 &0  \\
0& 0  & T\\
0 & U & 0  
\end{array}
\right ), \text{~~where $S= \left (\begin{array}{ll} 1 & 1\\
1& 0
\end{array}
\right)$, and}\\
T=\mathrm{diag} \{2,2,2,1,1,1,1,1,1,2,1 \},\, U=\mathrm{diag} \{2,2,1,1,2,2,2,2,2,1,2 \}. 
\end{array}
\]

There is a~$3|6$-structure on $\fg(2,3)^{(1)}/\fz $ that we express in the basis (\ref{basg23}) as follows:
\[
h_2^{[3]}=h_2, \, h_3^{[3]}=h_3,\, \quad  v^{[3]}=w^{[6]}=0 \text{~~for any $v$ (even) and $w$ (odd) root vectors in (\ref{basg23})}.
\]

\parbegin{Claim} \textup{(\cite{BGLL1})} The  space $\mathrm{H}^1 (\fg(2,3)^{(1)}/\fz; \fg(2,3)^{(1)}/\fz)$ is spanned by the $7$ even  cocycles of which we indicate the ones we need: 
\begin{equation} \label{g23cocycles}
\tiny
\begin{array}{ll}
\deg=-3:&\mathscr{D}_{-3}^1= 2\
 x_3\otimes \widehat x_8+x_6\otimes \widehat x_{10}+y_1\otimes \widehat x_4+
 y_4\otimes \widehat x_1+2\ y_8\otimes \widehat y_3+y_{10}\otimes \widehat y_6,\\[2mm]


\deg=0:&\mathscr{D}_{0}^1=x_2\otimes \widehat x_1+2\ x_6\otimes \widehat x_5+2\ x_9\otimes d
x_8+2\ y_1\otimes \widehat y_2+y_5\otimes \widehat y_6+y_8\otimes \widehat y_9,\\[2mm]

&\mathscr{D}_{0}^2=x_1\otimes \widehat x_1+2\ x_2\otimes \widehat x_2+2\ x_5\otimes d
x_5+x_6\otimes \widehat x_6+2\
x_8\otimes \widehat x_8+x_9\otimes \widehat x_9+\\

&2\ y_1\otimes \widehat y_1+y_2\otimes \widehat y_2+y_5\otimes \widehat y_5+2\ y_6\otimes
\widehat y_6+y_8\otimes \widehat y_8+2\
y_9\otimes \widehat y_9,\\


   

 \end{array}
 \end{equation}
 \end{Claim}
We easily see that $ \mathfrak{out}(\fg(2,3)^{(1)}/\fz)\simeq \mathfrak{psl}(3)$. 

A direct computation shows that the bilinear form $\mathscr{B}_\fa$ is  $\mathscr{D}$-invariant for the derivations  $\mathscr{D}$ given in (\ref{g23cocycles}).
The isomorphism  
\[
x_1\leftrightarrow y_1, \quad x_2\leftrightarrow y_2, \quad x_3\rightarrow y_3, \quad y_3\rightarrow 2x_3,
\]
permutes the derivations $\mathscr{D}^{1}_{-3}\leftrightarrow \mathscr{D}^{1}_{3}$ introducing $\mathscr{D}^{1}_{3}$. The isomorphism 
\begin{equation}
x_1\leftrightarrow x_2, \quad y_1\leftrightarrow y_2, \quad x_3\leftrightarrow x_3,  \quad y_3\leftrightarrow y_3,
\end{equation}
permutes the derivations $\mathscr{D}^{1}_{-3} \leftrightarrow\mathscr{D}^{2}_{-3}$, $\mathscr{D}^{1}_{0} \leftrightarrow \mathscr{D}^{3}_{0}$, and $\mathscr{D}^{1}_{3}\leftrightarrow \mathscr{D}^{2}_{3}$ turning the ones defined in \eqref{g23cocycles} and $\mathscr{D}^{1}_{3}$ into new ones. The derivations $\mathscr{D}^{1}_{-3}, \mathscr{D}^{1}_{0}$  and $ \mathscr{D}^{2}_{0}$ have the $p$-property.

Let us describe now the cubic forms. We fix an ordered basis of $\fa_\ev$ as follows
\[
e_1=h_2, e_2=h_3, e_3=x_1, e_4=x_2, e_5=x_4, e_6= x_{11},  e_7=y_1, e_8=y_2, e_9=y_4, e_{10}= y_{11}.
\]
For every $a=\sum_{1\leq i \leq 10}\lambda_i e_i$, we define 
\begin{eqnarray}
\label{LP1}\mathscr{P}_1(a)&:=&  \lambda_1 \lambda_3 \lambda_5+ \lambda_2 \lambda_3\lambda_5+ \lambda_3^2\lambda_4+ 2 \lambda_5^2 \lambda_8, \\[1mm]
\label{LP2}\mathscr{P}_2(a)&:=& 2\lambda_1\lambda_3\lambda_8+ 2 \lambda_2 \lambda_3 \lambda_8+ 2 \lambda_3^2 \lambda_9+ 2 \lambda_5 \lambda_8^2, \\[1mm]
\label{LP3}\mathscr{P}_3(a)&:=& 2 \lambda_1 \lambda_3 \lambda_7+ \lambda_1 \lambda_4 \lambda_8+ 2 \lambda_2 \lambda_3\lambda_7+ \lambda_2 \lambda_4\lambda_8+ 2 \lambda_3 \lambda_4 \lambda_9+ \lambda_5\lambda_7\lambda_8.
\end{eqnarray}
For each of these cubic forms, a~direct computation shows that the condition (\ref{Qua}) is satisfied. 
Let us summarize and introduce notation:
\begin{equation}\label{g23}
\footnotesize
\renewcommand{\arraystretch}{1.4}
\begin{tabular}{|c|c|c|c|c|} \hline
Derivation & $\mathscr{P}(a)$&$\gamma$ &$a_0$&
Double extension\\
\hline
$\mathscr{D}_{-3}^1$ & see Eq. (\ref{LP1}) &0&0& $\widetilde \fg(2,3)$ \\ \hline
   $\mathscr{D}_0^1$ & see Eq. (\ref{LP2})&0&0& $ \widehat{\fg}(2,3)$ \\ \hline
$\mathscr{D}_0^2$ & see Eq. (\ref{LP3}) &1& 0&$\fg(2,3)$\\ \hline
\end{tabular}
\end{equation}

\parbegin{Claim}  $\dim \mathrm{H}^3(\fg(2,3))=15$, $\dim \mathrm{H}^3(\widetilde \fg(2,3))=20$ and $\dim \mathrm{H}^3(\widehat{\fg}(2,3))=18$, hence  $\fg(2,3), \widetilde{\fg}(2,3)$ and $\widehat{\fg}(2,3)$ are pairwise not isomorphic. 
\end{Claim}
\subsubsection{$\fg(3,3)^{(1)}/\fz$, see \cite{BGL}}
Consider the Lie superalgebra  $\fg(3,3)$ with the Cartan matrix
\begin{equation*}\label{eq2.4.2b}\tiny
\begin{pmatrix}
 2 & -1 & 0 & 0 \\
 -1 & 2 & -1 & 0 \\
 0 & -2 & 2 & -1 \\
 0 & 0 & 1 & 0\end{pmatrix}\ \text{ with positive roots }\tiny\begin{array}{l}
 x_1,\ x_2,\ x_3,\ x_4,\\
x_5=[x_1, x_2],\ x_6=[x_2, x_3],\ x_7=[x_3, x_4],\\
x_8=[x_3, [x_1, x_2]],\ x_9=[x_3, [x_2, x_3]],\ x_{10}=[x_4, [x_2, x_3]],\\
x_{11}=[x_3, [x_3, [x_1, x_2]]], \ x_{12}=[[x_1, x_2], [x_3, x_4]],
x_{13}=[[x_2,
x_3], [x_3, x_4]], \\
x_{14}=[[x_2, x_3], [x_3, [x_1, x_2]]], \ \ x_{15}=[[x_3, x_4], [x_3,
[x_1, x_2]]], \\
x_{16}=[[x_3, [x_1, x_2]], [x_4, [x_2, x_3]]], \\
x_{17}=[[x_4, [x_2, x_3]], [x_3, [x_3, [x_1, x_2]]]].\end{array}
\end{equation*} 
Denote by $y$'s the corresponding negative roots, and set $h_i = [x_i,y_i]$ for $i=1,2,3$.

The simple core $\fa:=\fg(3,3)^{(1)}/\fz $ admits a~NIS given in the ordered basis 
\begin{equation}
\label{basg33}
h_2:=[x_2,y_2], \; h_3:= [x_3,y_3],\;  h_4:= [x_4,y_4],\;  x_1,\ldots x_{17},\; y_1, \ldots,y_{17},
\end{equation}
by the Gram matrix
\[
\mathscr{B}_{\fa}=\left (
\begin{array}{ccc}
S& 0 &0  \\
0& 0  & T\\
0 & U & 0  
\end{array}
\right ), \text{~~where}
\]
\[
\begin{array}{lcl}
T&=&\mathrm{diag} \{2,2,1,2,1,2,1,1,2,2,1,1,1,2,2,1,2, \},\\[2mm]
U&=&\mathrm{diag} \{2,2,1,1,1,2,2,1,2,1,1,2,2,2,1,2,1\},\\[2mm]
S&=&E^{1,1}+2 E^{2,1}+2 E^{2,1}+2 E^{2,2}+E^{3,1}+E^{1,3}.
\end{array}
\]

There is a~$3|6$-structure on $\fg(2,3)^{(1)}/\fz $ that we express in the basis (\ref{basg33}) as follows:
\[
h_2^{[3]}=h_2, \, h_3^{[3]}=h_3,\, \, h_4^{[3]}=h_4,\, \quad  v^{[3]}=w^{[6]}=0 \text{ if $v$ is even and $w$ is odd roots in (\ref{basg33})}.
\]

\parbegin{Claim} \textup{(\cite{BGLL1})} The  space $\mathrm{H}^1(\fg(3,3)^{(1)}/\fz;\fg(3,3)^{(1)}/\fz)$is spanned by the following cocycles:
\[\tiny
\begin{array}{ll}
\deg=-8:&\mathscr{D}_{-8}= y_4\otimes \widehat x_{17} + y_7\otimes \widehat x_{16} + 2\
y_{10}\otimes \widehat x_{15} +y_{12}\otimes \widehat x_{13} +y_{13}\otimes \widehat 
x_{12} + 2\ y_{15}\otimes \widehat 
x_{10}+y_{16}\otimes \widehat x_7 +y_{17}\otimes \widehat x_4, \\[2mm]

\deg=0:&\mathscr{D}_{0}=2\ x_4\otimes \widehat x_4 + 2\ x_7\otimes \widehat x_7 + 2\ x_{10}\otimes \widehat 
x_{10} + 2\ x_{12}\otimes \widehat x_{12} + 2\
x_{13}\otimes \widehat x_{13} + 2\ x_{15}\otimes \widehat x_{15} + 2\ x_{16}\otimes \widehat x_{16} \\

&+ 2\ x_{17}\otimes \widehat x_{17} + y_4\otimes
\widehat y_4 +y_7\otimes \widehat y_7 +y_{10}\otimes \widehat y_{10}
+y_{12}\otimes \widehat y_{12} +y_{13}\otimes \widehat y_{13} + y_{15}\otimes \widehat y_{15} +y_{16}\otimes \widehat y_{16} +y_{17}\otimes \widehat 
y_{17},\\[2mm]
\deg=8:&\mathscr{D}_{8}= x_4\otimes \widehat{y}_{17} +
x_7\otimes \widehat{y}_{16} 
+ 2\,
   x_{10}\otimes \widehat{y}_{15}
  +
   x_{12}\otimes \widehat{y}_{13}
   +
 x_{13} \otimes \widehat{y}_{12}  + 
 2\,
   x_{15}\otimes \widehat{y}_{10}
   +
   x_{16} \otimes \widehat{y}_7
   +
 x_{17} \otimes \widehat{y}_4 .
\end{array}
\]
\end{Claim}
We will omit the details here. By symmetry, the cases $\mathscr{D}_{-8}$ and $\mathscr{D}_{8}$ give isomorphic double extensions. There are two non-isomorphic double extensions of $\fg(3,3)^{(1)}/\fz$ given as in the table below:
\begin{equation}\label{g33}
\footnotesize
\renewcommand{\arraystretch}{1.4}
\begin{tabular}{|c|c|c|c|c|c|} \hline
Derivation & The cubic form $\mathscr{P}$&$\gamma$ &$a_0$&
Double extension (DE)& $ \mathrm{dim } \, \mathrm H^2(\text{DE})$\\
\hline
$\mathscr{D}_{-8}$ & $0$ &0&0& $\widetilde \fg(3,3)$ &1 \\ \hline
   $\mathscr{D}_0$ & $0$ &0&0& $ \fg(3,3)$  & 0\\ \hline
\end{tabular}
\end{equation}
\subsection{$\mathfrak{psq}(n)$ for $n>2$ and $p\neq 2$} The case $p=2$ is completely different: there is a $1|1$-dimensional space of NISes on $\mathfrak{psq}(n)$, see \cite{KLLS}; to classify double extensions of $\mathfrak{psq}(n)$ is an open problem; to solve it, one should extend the results of \cite{BGLL1}. For completeness, note that for $p=0$, the double extension of $\mathfrak{psq}(n)$ is $\fq(n)$.

\sssbegin{Conjecture} \textup{(Verified for $n=3,4$ and $5$ and $p=3,5$)} The  space $\mathrm{H}^1 (\mathfrak{psq}(n);  \mathfrak{psq}(n))$ is spanned by one odd cocycle. Therefore, the only double extension of $\mathfrak{psq}(n)$ is $\mathfrak{q}(n)$.
\end{Conjecture}
\section{Appendix: Pre-Lie superalgebras in characteristic 3}\label{sec: appen}
A Lie superalgebra $\fa$ over a~field of characteristic $p$, where $p\not =2,3$, is a~$\Bbb Z/2 \Bbb Z $-graded vector space $\fa=\fa_\ev\oplus \fa_\od$, endowed with a~bilinear operation $[-, -]$ that satisfies anti-commutativity and Jacoby identities amended by the Sign Rule.

 The~{\bf Lie superalgbera $\fa$ in characterstic 3} must satisfy one more condition (\cite{LeP})
\begin{equation}
\label{JIp=3odd} [a, [a,a]]=0 \quad \text{for all $a\in \fa_\od$}.
\end{equation}
Obviously, condition (\ref{JIp=3odd}) is a~consequence of the Jacobi identity  if $p\not =3$. Accordingly, the 2-cocycle, say $\omega$, 
describing central extensions of $\fa$ should satisfy the following conditions:
\begin{eqnarray}
\label{Ext1}\omega(a,b)=-(-1)^{p(a)p(b)} \omega (b,a) \quad \text{for all $a,b\in \fa$},\\[2mm]
\label{Ext2}\omega(a, [b,c])= \omega([a,b],c)+(-1)^{p(a)p(b)} \omega(b, [a,c]) \quad \text{for all $a,b,c\in \fa$,}\\[2mm]
\label{Ext3} \text{(Only for $p=3$)} \quad  \omega (a, [a,a])=0 \quad \text{for all $a\in \fa_\od$}.
\end{eqnarray}
It was shown in \cite{BGLL1} that 
\[
\mathrm{H}^2(\fosp(1|2))=0  \quad \text{and} \quad \mathrm{H}^2(\fbrj(2;3))=0.
\] 
Since $\fosp(1|2)$ and $\fbrj(2;3)$ do not admit non-trivial central extensions, they have no double extentions. 

On the other hand, in Subsections \ref{prebrj} and \ref{preosp}, we observe the following interesting fact: both $\fosp(1|2)$ and $\fbrj(2;3)$ are NISs superalgebras and admit an
odd outer derivation $\mathscr{D}$, given explicitly (recall that $p=3$) in (\ref{cocybrj3}) and (\ref{cocyposp}). For any of these $\mathscr{D}$, the map
\[
\omega(a,b):=\mathscr{B}(\mathscr{D}(a),b)
\]
does satisfy the 2-cocycle conditions (\ref{Ext1}) and (\ref{Ext2}), but not the condition (\ref{Ext3}) since for any of these $\mathscr{D}$ we have ${\mathscr{B}(\mathscr{D}(a),[a,a])\not =0}$ for some $a\in \fa_\od$. Moreover, there is no $\alpha\in \fa^*$ such that $\omega(a,b)=\alpha([a,b]);$ otherwise, this would have implied that $\mathscr{D}$ was an inner derivation.

It follows that the bracket, see Eq. (\ref{*}), defined on the double extensions of $\fosp(1|2)$ and $\fbrj(2;3)$ satisfies the Jacoby identity but not (\ref{JIp=3odd}). Leites suggested to call these algebras for which condition (\ref{JIp=3odd}) is violated, but the Jacoby identity holds, \textit{pre-Lie superalgebras}. 

The double extensions of $\fosp(1|2)$ and $\fbrj(2;3)$, both for $p=3$, are \textbf{pre-Lie} superalgebras. Currently, these are the {\bf only} examples known in the literature of central extensions of simple Lie superalgebras which are pre-Lie superalgebras but not Lie superalgebras. 

A similar situation occurs in characteristic 2, where a~central extension of $\fh_I(n; \underline{N})$ is not a~Lie algebra, but a~\textbf{Leibniz} algebra, see \cite{BGLLS}.


\section*{Acknowledgements} We heartily thank A.~Lebedev for clarifying to us that the double extension of $\fosp(1|2)$ in 
characteristic 3 is not a~Lie superalgebra, but a~pre-Lie superalgebra, and the referee for meticulous job and huge help. We also thank D.~Leites for discussions 
and clarifications.  S.~Bouarroudj was supported by the grant AD 065 NYUAD. 

\end{document}